\documentclass{article}
\usepackage{amsfonts, amsmath, amssymb, amscd,amsthm, graphicx,latexsym,epsfig} 

\usepackage[all]{xy}

\usepackage{auto-pst-pdf}

\newtheorem{thm}{Theorem}[section]

\newtheorem{cry}{Corollary}[section]
\newtheorem{lem}{Lemma}

\newtheorem{prop}{Proposition}[section]
\newtheorem{defin}{Definition}[section]

\newcommand{\fg}[1]{\pi_1(\projective\moins #1 )}

\newcommand{\onto}{\twoheadrightarrow}

\newcommand{\mon}{^{-1}}
\newcommand{\mtw}{^{-2}}
\newcommand{\mth}{^{-3}}

\newcommand{\inv}{^{-1}}
\newcommand{\rel}[1]{{\cal R}_{#1}}
\newcommand{\complex}{\mathbb{C}}
\newcommand{\projective}{\mathbb{P}^2}
\newcommand{\sphere}{\mathbb{P}^1}

\newcommand{\realline}{\mathbb{R}}

\newcommand{\naturals}{\mathbb{N}}
\newcommand{\integers}{\mathbb{Z}}
\newcommand{\braids}{\mathbb{B}}
\newcommand{\moins}{\backslash}

\newcommand{\ffigure}[1]{(Figure~\ref{#1})}
\newcommand{\fffigure}[1]{Figure~(\ref{#1})}

\makeatletter
\long\def\@makecaption#1#2{%
\vskip 10\p@
\setbox\@tempboxa\hbox{#1\ \  #2}%
\ifdim \wd\@tempboxa >\hsize
#1\ \#2\par
\else
\hbox to \hsize{\hfil\box \@tempboxa\hfil}%
\fi}
\makeatother

\newcommand{\figur}[1]{
\begin{figure}\label{#1}
\begin{center}
\includegraphics{#1.ps}
\caption{\label{#1}}
\end{center} 
\end{figure}}

\begin{document}

\markboth{A. Muhammed Uluda\u g}
{Fundamental Groups of a Class of Rational Cuspidal Plane Curves}

\title{Fundamental Groups of a Class of Rational Cuspidal Plane Curves}

\author{A. Muhammed Uluda\u{g}\footnote{
{Galatasaray University, Department of Mathematics,}
{\c{C}{\i}ra\u{g}an Cad. No. 36, 34357 Be\c{s}ikta\c{s}}
{\.{I}stanbul, Turkey}}}

\maketitle

\begin{abstract}
We compute the presentations of fundamental groups of the complements of a class of rational cuspidal projective plane curves classified by Flenner, Zaidenberg, Fenske and Saito. We use the Zariski-Van Kampen algorithm and exploit the Cremona transformations used in the construction of these curves. We simplify and study these group presentations so obtained and determine if they are abelian, finite or big, i.e. if they contain free non-abelian subgroups. We also study the quotients of these groups to some extend.
\end{abstract}


\section{Introduction}
A projective plane curve $C\subset \projective$ 
is said to be \textit{cuspidal} if all its
singularities are irreducible. 
It is said to be of \textit{type} $(d,m)$ if $C$ is of degree
$d$ and the maximal multiplicity of its singularities is 
$m$, i.e. $m=\max_{p\in  C}(\mbox{mult}_p C)$. 
Rational cuspidal curves of type $(d,d-2)$ and those of type $(d,d-3)$, having
at least three cusps has been
classified by Flenner and Zaidenberg. 
It turns out that these curves have exactly three cusps, moreover, 
the rational 3-cuspidal 
curves of type $(d,d-2)$ can be obtained from  the quadric 
$\{xy-z^2=0\}$, and the rational 3-cuspidal 
curves of type $(d,d-3)$ can be obtained from
the cubic $\{xy^2-z^3=0\}$ by means of Cremona transformations.
Our aim in this paper is to compute the fundamental groups of complements of these curves as well  as those ones constructed in the ensuing work of Fenske, Sakai and Tono. 

Given an algebraic curve $C\subset \projective$, it is interesting to
study  the topology of its complement $\projective\moins C$
from the point of view of the classification theory of algebraic curves. 
The idea is an analogue of the leading principle in the knot theory: in
order to understand a knot ${\cal K}\subset S^3$, look at the topology
of the knot complement $S^3\moins {\cal K}$.
Another strong motivation for studying 
$\fg{C}$ comes from the surface theory: 
A good knowledge of $\fg{C}$ allows one to construct  
Galois coverings of $\projective$ branched at $C$ (see \cite{uludagk3} for explicit examples).
Moreover, if $X\rightarrow \projective$ is a 
branched covering, with $C$ as the
branching locus, one can hope to derive some invariants of $X$ from
the invariants of the topology of $\projective\moins C$. 
The study of $\fg{C}$ has been initiated by  Zariski in the thirties \cite{Zariskiotpe}.
Although the prevailing convention in the literature is to call  
$\pi_1(\projective \moins C)$, by abuse of language,
 the fundamental group \textit{of} $C$, following 
Degtyarev~\cite{Degtyarev1}, 
we shall take the liberty to call it simply 
\textit{the group of $C$}.
Degtyarev classified the possible groups of curves of degree $d$ with a
singular point of multiplicity $d-2$ and calculated the groups of all 
quintics (irreducible or not). Hence, groups of curves of degree $\leq 5$ are known.
This is not the case starting with the sextics.
A special class of sextics, namely the rational cuspidal ones, 
has been classified by Fenske~\cite{Fe1}, see Theorem~\ref{fe3}. 
Their groups are given in Corollary~\ref{crysextics}.

Below we review these classification results followed by Artal's and our results on the fundamental groups. Since group computations uses the explicit construction of these curves, we felt obliged to give a detailed account of classification results.

\paragraph{Acknowledgements.} These results are from authors's thesis. Hereby I express my gratitude towards M. Zaidenberg who suggested the problem. 
I am indebted to Alex Dimca for informing that they are useful and referred, and for his encouragement to publish them. 
I am thankful to Hakan Ayral for his help with the graphics.

\section{Classifications and fundamental groups}
We  use the following conventions settled in~\cite{FlZa2}:
the multiplicity sequence of a cusp will be
called the \textit{type} of this cusp. 
Recall that if
$$
V_{n+1}\rightarrow V_n\rightarrow\cdots V_1\rightarrow V_0=\complex^2
$$
is a minimal resolution of an irreducible analytic curve singularity
germ $(C,0)\subset (\complex^2, 0)$, and $(C_i, P_i)$ denotes the
proper transform of $(C,0)$ in $V_i$, so that $(C_0, P_0)=(C,O)$, 
then the sequence $$[m^{(n+1)}, m^{(n)},\dots ,m^{(1)}, m^{(0)}]$$,
where $m^{(i)}:={\rm mult}_{P_i}C_i$, is called the 
\textit{multiplicity sequence} of $(C,0)$. 
Evidently, $m^{(i+1)}\leq m^{(i)}$, $m^{(n)}\geq 2$, and $m^{(n+1)}=1$.
The (sub)sequence
$$
\underbrace{k,k,k,\dots ,k}_{m \;\scriptsize{\mbox{times}}}
$$
will be abbreviated by $k_m$. For instance, $[kn, k_{n+m},k-1]$ is the
sequence 
$[kn,\underbrace{k,k,k,\dots ,k}_{n+m \;\scriptsize{\mbox{times}}},
k-1]$.
Also, the last term of the multiplicity sequence will be omitted. 
Under this notation, $[2]$ corresponds to a simple cusp and $[2_3]=[2,2,2]$ 
corresponds to a  ramphoid cusp. 
\begin{thm}\textbf{(Flenner-Zaidenberg~\cite{FlZa1})}
\label{fz1}
A rational cuspidal curve of type $(d,d-2)$ with at least three cusps
has exactly three cusps. 
For each $(d,n,m)\in \naturals^3$ such that 
$d\geq 4$, $n\geq m>0$ and $n+m=d-2$
there is exactly one (up to projective
equivalence) such curve $C$, whose  
 cusps are of types 
 $[d-2]$, $[2_n]$, $[2, 2_m]$. There are no other rational cuspidal 
curves of type $(d,d-2)$ with number of cusps $\geq 3$. 
\end{thm}
Note that an irreducible curve of type $(d,d-1)$ is rational, and has
an abelian group by a classical computation due to Zariski.
\begin{thm}\textbf{(Artal~\cite{Artal1})}
\label{Artal}
Let $C$ be as in Theorem~\ref{fz1}. 
Then the  group of  $C$ admits the presentation
$$
\fg{C}=
\langle a,b\,|\, (ba)^{d-1}=b^{d-2},\quad (ba)^k b=a(ba)^k\rangle,
$$
where $k\geq 0$ and $2k+1={\rm gcd}(2n+1,\,2m+1)$. Hence, the group
depends only on $(d,k)$. This group is abelian if and only if $k=0$,
finite of order $12$ if $(d,n,m)=(4,1,1)$, finite of order $840$ if
$(d,n,m)=(7,4,1)$, otherwise it is a big group~\footnote{Recall that
  we call a group \textit{big} if it has a non-abelian free subgroup.}.
\end{thm}
Note that the case $(d,n,m)=(4,1,1)$ corresponds to the three cuspidal
quartic, whose group had been calculated already 
by Zariski~\cite{Zariskiotpe}.
\par
Taking for example $C_1$ with $(d,n,m)=(13, 10, 1)$ and $C_2$ with 
$(d,n,m)=(13, 7, 4)$, one has the following result.
\begin{cry}\textbf{(Artal~\cite{Artal1})}
There exist infinitely many pairs $(C_1, C_2)$ of curves with
isomorphic (big) groups, but 
different homeomorphism types of the pairs 
$(\projective, C_1)$, $(\projective, C_2)$.
\end{cry}
\begin{thm}\textbf{(Flenner-Zaidenberg~\cite{FlZa2})}
\label{fz2}
A rational cuspidal curve of type $(d,d-3)$ with at least three cusps
has exactly three cusps. 
For each $d=2n+3$, $n\geq 1$, there is exactly one (up to projective
equivalence) such curve $C_n$, 
and the cusps of $C_n$ are of types $[2n,2_n]$, $[3_n]$, $[2]$.
There are no other rational cuspidal 
curves of type $(d,d-2)$ with number of cusps $\geq 3$.
\end{thm}
\begin{thm}
\label{groupssflza}
Let $C_n$ be as in Theorem~\ref{fz2}, where $d=2n+3$. 
Then the  group of $C_n$ admits the presentation
$$
\pi_1(\projective\moins C_n)=
\langle c,b\,|\,cbc=bcb,\quad b^n c^{n+2}=c^{n+2}b^n,\quad
(b^{-n}cb^2)^{n+1}c^{n^2}=1\rangle.
$$
This group is big when $n$ is odd and $\geq 7$, abelian when
$n\in\{0,1,2\}$, and finite of order $8640$ for $n=3$, and finite of
order 1560 for $n=5$. 
\end{thm}
This is proved in 3.3.1 below.

These classification results have been completed (independently) in the works of Fenske and 
Sakai-Tono:
\begin{thm}\textbf{(Fenske~\cite{Fe1}, Sakai-Tono~\cite{SaTo})}
\label{fe1}
Let $n,m\in \naturals $ be such that $n\geq 1$ and $ 0\leq m <n$.
\\(i)
A rational cuspidal curve of type $(d,d-2)$ 
with exactly one cusp exists if and only if
$d=2n+2$. 
Such a curve is unique up to projective equivalence, 
and the type of its cusp is $[2n, 2_{2n}]$.
\\(ii)
A rational cuspidal curve of type $(d,d-2)$ 
with exactly two cusp exists if and only if
\\(a) $d=2n+3$, with types of cusps $[2n+1,2_n]$, $[2_{n+1}]$;
\\(b) $d=n+2$,  with types of cusps $[n]$, $[2_n]$;
\\(c) $d=2n+2$, with types of cusps $[2n, 2_{n+m}]$, $[2_{n-m}]$.\\
Moreover, all these curves are unique up to projective equivalence.
\end{thm}
\begin{thm}\textbf{(Fenske~\cite{Fe1})}
\label{fe2}
\\Let $n,m\in \naturals $ be such that $n\geq 1$ and $ 0\leq m <n$.
\\(i)
A rational cuspidal curve of type $(d,d-3)$ 
with exactly one cusp exists if and only if
\\(a) $d=3n+3$, where the type of the cusp is $[3n, 3_{2n},2]$;
\\(b) $d=5$,    where the type of the cusp is $[2_6]$.
These curves are unique up to projective equivalence.
\\(ii) The only existing rational cuspidal curves of type $(d,d-3)$
with exactly two cusps are the following ones:
\begin{center}
\begin{tabular}{|l|r|r|}
\hline
      &  degree &  types of cusps               \\ \hline
$1$   &  $7$    &  $[4]\,,\,[3_3]$                    \\ \hline
$2$   &  $6$    &  $[3]\,,\, [3_2,2]$                 \\ \hline    
$3$   &  $5$    &  $[2_4]\,,\, [2_2]$                 \\ \hline
$4$   & $2n+3$  &  $[2n,2_n]\,,\, [3_n,2]$            \\ \hline
$5$   & $2n+4$  &  $[2n+1,2_n]\,,\, [3_{n+1}]$        \\ \hline    
$6$   & $2n+3$  &  $[2n,2_{n+1}]\,,\, [3_n]$          \\ \hline    
$7$   & $3n+3$  &  $[3n,3_{2n}]\,,\, [2]$             \\ \hline    
$8$   & $3n+4$  &  $[3n+1,3_n]\,,\, [3_{n+1}]$        \\ \hline
$9$   & $3n+3$  &  $[3n,3_{n+m},2]\,,\, [3_{n-m}]$    \\ \hline
${10}$& $3n+3$  &  $[3n,3_{n+m}]\,,\, [3_{n-m},2]$    \\ \hline  
${11}$& $3n+3$  &  $[3n+2,3_n,2]\,,\, [3_{n+1}, 2]$   \\ \hline
\end{tabular}   
\end{center}
These curves are unique up to projective equivalence.
\end{thm}
Considering the cases $d=6$ in the above theorems leads to 
a complete  classification of rational cuspidal curves of degree 6,
given in~\cite{Fe1}. 
\begin{thm}\label{fe3} \textbf{(Fenske \cite{Fe1})}
Up to projective equivalence, 
rational cuspidal curves of degree 6 are the following ones:
\begin{center}
\begin{tabular}{|l|r|r|}
\hline
        &    types of cusps               \\ \hline
$1$   &  $[5]$                    \\ \hline
$2$   &  $[4,2_4]$                 \\ \hline    
$3$   &  $[3_3, 2]$                 \\ \hline
$4$   & $[3_3]\,,\, [2]$            \\ \hline
$5$   & $[3_2,2]\,,\, [3]$        \\ \hline    
$6$   & $[3_2]\,,\, [3,2]$          \\ \hline    
\end{tabular}   \quad
\begin{tabular}{|l|r|r|}
\hline
        &    types of cusps               \\ \hline
$7$   & $[4,2_3]\,,\, [2]$             \\ \hline    
$8$   & $[4,2_2]\,,\, [2_2]$        \\ \hline
$9$   & $[4]\,,\, [2_4]$    \\ \hline
${10}$& $[4]\,,\, [2_3], [2]$    \\ \hline  
${11}$& $[4]\,,\, [2_2], [2_2]$   \\ \hline
&\\ \hline
\end{tabular}  
\end{center}

\end{thm}
\par
In fact, in \cite{Fe1} the following 
longerl list of curves is constructed.
\begin{thm} \label{fe4}
\textbf{(Fenske~\cite{Fe1})} Let $d\geq 2$ and $0\leq m<n$ be
  integers. The following rational cuspidal plane curves exist:
 \begin{center}
\begin{tabular}{|l|r|r|}
\hline
      &type of curves  $(d,m)$      &  type of cusps\\ \hline
$1$     &  $(kn+k,kn)$       &  $[kn, k_{n+m},k-1]\,,\,[k_{n-m}]$\\ \hline
$1a$    &  $(kn+k,kn)$       &  $[kn,k_{2n},k-1] $\\ \hline    
$2$     &  $(kn+k,kn)$       &  $[kn,k_{n+m}]\,,\,[k_{n-m},k-1]$\\ \hline
$2a$    &  $(kn+k,kn)$       &  $[kn,k_{2n}]\,,\,[k-1]$\\ \hline
$3$     &  $(kn+k+1,kn+1)$   &  $[kn+1,k_{n}]\,,\,[k_{n+1}]$\\ \hline    
$4$     &  $(kn+k+1,kn)$     &  $[kn,k_{n+1}]\,,\,[(k+1)_n]$\\ \hline    
$5$     &  $(kn+k+1,kn)$     &  $[kn,k_n]\,,\,[(k+1)_n,k]$\\ \hline    
$6$     &  $(kn+k+2,kn+1)$   &  $[kn+1,k_n]\,,\,[(k+1)_{n+1}]$\\ \hline
$7$     &  $(kn+2k-1,kn+k-1)$&  $[kn+k-1,k_n,k-1]\,,\,[k_{n+1},k-1]$\\ \hline
$8$     &  $(n+2,n)$         &  $[n]\,,\,[2_n]$\\ \hline  
 \end{tabular}   
\end{center}
\end{thm}
Beware of the following exception in the table above: In case $n=1$,
the curves (4) and (5) are of type $(kn+k+1, k+1)$, instead of 
$(kn+k+1, nk)$.

\begin{thm}
\label{table}
Groups of the curves in Theorem~\ref{fe4} are as follows:
\begin{center}
\begin{tabular}{|l|r|r|}
\hline
      &       \\ \hline
$1$     &  $\langle\alpha, \beta,
y\;|\;\alpha=y^{-k}\beta^{k-1},\quad[\beta, y^k]=\alpha(\alpha\beta)^m=
\beta(\alpha\beta)^{n-m}=1\rangle$                 \\ \hline
$1a$    &  abelian                           \\ \hline    
$2$     &  $\langle\alpha, \beta, y\;|\;
\alpha=y^{-k}\beta^{k-1},\quad
[\beta, y^k]=  
\alpha(\alpha\beta)^{n-m}=
\beta(\alpha\beta)^{m}=1\rangle$                   \\ \hline
$2a$    &  $\langle\beta,
y\;|\;y^k=\beta^{k-1},\quad\beta^{n+1}=1\rangle$        \\ \hline
$3$     &   abelian                          \\ \hline    
$4$     &   $\langle x,y\,|\, 
(xy)^{n+1}y^{n^2}=[y, (xy)^k]=1,\quad x^ky^nx=(xy)^k\rangle$
\\ \hline    
$5$     &   abelian                          \\ \hline    
$6$     &   abelian                          \\ \hline
$7$     &   abelian                          \\ \hline
$8$     &   abelian                          \\ \hline  
 \end{tabular}   
\end{center}
Groups (1) are central extensions of the group $\integers_k
*\integers_j$, where  $j:=\mbox{g.c.d.}(mk+k-1,n+1)$. 
Thus, they are abelian if $j=1$, 
and big if $j\geq 2$. Groups (2) are central extensions of the group $\integers_k
*\integers_j$, where $j:=\mbox{g.c.d.}(1+mk, n+1)$. 
Thus, they are abelian if $j=1$, and big if $j\geq 2$.
The same conclusion is true for the groups (2a), where this time
$j:=\mbox{g.c.d.}(k-1, n+1)$.
Groups (4) are abelian if $j:=(n+1, k)=1$ or $n=1$. 
Otherwise, they are big with the following exceptions:
\\(i) If $(n,k)=(3,2)$, then the group is finite non-abelian of
order $72$ (the degree of the curve is $9$).
\\(ii) If $(n,k)=(5,2)$, then the group is finite non-abelian of
order $1560$ (the degree of the curve is $13$)
\\(iii) If (n,k)=(2,3), then the group is finite non-abelian of
order $240$ (the degree of the curve is $10$).
\end{thm}
This is proved in 3.3.2-3.3.8 below.
\begin{cry}  \label{cryd-2}
(i)  The group of a rational unicuspidal curve of
  type $(d,d-2)$ is abelian, 
\\
(ii) The group of a rational two-cuspidal curve $C$ of type 
$(d,d-2)$ is abelian unless $C$ is one of the curves described in 
Theorem~\ref{fe1}-(ii)c, and $j:=\mbox{g.c.d.}(2m+1, n+1)\neq 1$. 
In this case, the  group of $C$ is the big group with
the following presentation 
$$
\langle y,\beta\,|\,
[\beta, y^2]=y^{-2m-2}\beta^{2m+1}=y^{2m-2n}\beta^{2n-2m+1}=1\rangle.
$$
This group is a central extension of $\integers_2*\integers_j$.
\end{cry}
\proof

(i) This is the case (1a) with $k=2$ in Theorem~\ref{table}.\\
(ii) (a)  This is the case (3) with $k=2$ in Theorem~\ref{table}.
\\
(b) This is the case (8)  in Theorem~\ref{table}.
\\
(c)  This is the case (1) with $k=2$ in Theorem~\ref{table}. One
obtains the presentation easily by substituting $\alpha=y\mth \beta$.

\begin{cry} \label{cryd-3}
(i)  The group of a rational unicuspidal curve of type 
$(d,d-3)$  is abelian.
\\
(ii)  Groups of  rational two-cuspidal curves of type 
$(d,d-3)$ are given below:
\begin{center}
\begin{tabular}{|l|r|r|}
\hline
      &  degree &  group               \\ \hline
$1$   &  $7$    &  abelian                    \\ \hline
$2$   &  $6$    &  $\integers_2*\integers_3$                 \\ \hline    
$3$   &  $5$    &  abelian                 \\ \hline
$4$   & $2n+3$  &  abelian            \\ \hline
$5$   & $2n+4$  &  abelian        \\ \hline    
$6$   & $2n+3$  &  $\langle x,y\,|\, xy^nx=yxy,\quad
(xy)^{n+1}y^{n^2}=[y,xyx]=1\rangle$          \\ \hline    
$7$   & $3n+3$  &  $\langle y, \beta\,|\, \beta^2=y^3, \quad
\beta^{n+1}=1\rangle$             \\ \hline    
$8$   & $3n+4$  &  abelian        \\ \hline
$9$   & $3n+3$  &  $\langle\alpha, \beta,
y\;|\;\alpha=y^{-3}\beta^{2},\quad[\beta, y^3]=\alpha(\alpha\beta)^m=
\beta(\alpha\beta)^{n-m}=1\rangle$    \\ \hline
${10}$& $3n+3$  &  $\langle\alpha, \beta, y\;|\;
\alpha=y^{-3}\beta^{2},\quad
[\beta, y^3]=  
\alpha(\alpha\beta)^{n-m}=
\beta(\alpha\beta)^{m}=1\rangle$    \\ \hline  
${11}$& $3n+5$  &  abelian   \\ \hline
\end{tabular}
\end{center}
Groups  (6) are abelian if  $n$ is even or $n=1$. 
Otherwise they are big unless $n=3$ or $n=5$, in these cases the group
is finite of order $72$ and $1560$ respectively.
Groups (7) are abelian if $n$ is even, and big otherwise.
Groups (9) are abelian if $g.c.d.\,(3m+2, n+1)=1$, and big
otherwise.
Groups (10) are abelian if $g.c.d.\,(3m+1, n+1)=1$, and big otherwise.

\end{cry}
\proof
\
1. This is the case (4) in Theorem~\ref{table} with $k=3$ an
$n=1$. Hence, the group has the presentation
$$
\langle x,y\, |\, x^3yx=(xy)^3,\quad (xy)^2y=[y,(xy)^3]=1\rangle,
$$
which is easily seen to be abelian, by substituting $(xy)^2=y\mon$ in
the commutation relation.
\\
2. This is the case (1) in Theorem~\ref{table} with $k=3$, $n=1$, and 
$m=0$. The group is
$$ 
\langle \alpha \beta,y\, |\,\alpha=y\mth \beta^2, \quad
[\beta, y^3]=\alpha=\beta\alpha\beta=1 \rangle.
$$
Thus, $y^3=\beta^2=1$, i.e. the group is $\integers_2*\integers_3$.
\\
3. The group of this curve was found to be abelian by
Degtyarev~\cite{Degtyarev1}. 
\\
4. This is the case (5) with $k=2$ in Theorem~\ref{table}. The group
is thus abelian.
\\
5. This is the case (6) with $k=2$ in Theorem~\ref{table}, and the
group is abelian.
\\
6. This is the case (4) with $k=2$ in Theorem~\ref{table}. 
\\
7. This is the case (2a) with $k=3$ in Theorem~\ref{table}.
\\
8. This is the case (3) with $k=3$ in Theorem~\ref{table}.
\\
9. This is the case (1) with $k=3$ in Theorem~\ref{table}.
\\
10. This is the case (2) with $k=3$ in Theorem~\ref{table}.
\\
11. This is the case (1) with $k=3$ in Theorem~\ref{table}.

\begin{cry}  
\label{crysextics}
Groups of rational cuspidal sextics are listed below.
\begin{center}
\begin{tabular}{|l|r|r|}
\hline
        &    types of cusps  &   group          \\ \hline
$1$   &  $[5]$ & abelian                    \\ \hline
$2$   &  $[4,2_4]$ &abelian                 \\ \hline    
$3$   &  $[3_3, 2]$ &abelian                 \\ \hline
$4$   &  $[3_3], [2]$ & $\integers_2*\integers_3$            \\ \hline
$5$   & $[3_2,2]\,,\, [3]$ & $\integers_2*\integers_3$       \\ \hline    
$6$   & $[3_2]\,,\, [3,2]$ & $\integers_2*\integers_3$         \\ \hline    
$7$   & $[4,2_3]\,,\, [2]$ & $\integers_2*\integers_3$            \\ \hline    
$8$   & $[4,2_2]\,,\, [2_2]$ & abelian        \\ \hline
$9$   & $[4]\,,\, [2_4]$ & abelian    \\ \hline
${10}$& $[4]\,,\, [2_3], [2]$  & abelian   \\ \hline  
${11}$& $[4]\,,\, [2_2], [2_2]$ & $\langle a,b\,|\, (ab)^5=b^4,\quad
a(ba)^2=(ba)^2 b $  \\ \hline
\end{tabular}   
\end{center}
\end{cry}
\proof

1-2-3: These sextics are unicuspidal, hence their groups are abelian
by Corollaries \ref{cryd-2} and \ref{cryd-3}.
\\
4. This is the case (7) in Corollary \ref{cryd-3} with $n=1$.
\\
5. This is the case (2) in  Corollary \ref{cryd-3}.
\\
6. This is the case (10) in  Corollary \ref{cryd-3} with $n=1$ and $m=0$. 
Thus $j=1$, and the group is abelian.
\\
7. This is the curve in  Corollary \ref{cryd-2}-(c) with $n=2$ and $m=1$. 
Thus $j=3$, and the group has the presentation
$$
\langle y,\beta\,|\,
[\beta, y^2]=y^{-4}\beta^{3}=y^{-2}\beta^{3}=1\rangle,
$$
which easily seen to be isomorphic to $\integers_2*\integers_3$.
\\
8. This is the curve in  Corollary \ref{cryd-2}-(c) with $n=2$ and $m=0$. 
Thus $j=1$, and the group is abelian.
\\
9. This is the curve in Corollary \ref{cryd-2}-(b).
\\
10. This is the curve in Theorem~\ref{Artal} with
$(d,n,m)=(6,3,1)$. Thus, $k=1$, and the group is abelian.
\\
11.
This is the curve in Theorem~\ref{Artal} with
$(d,n,m)=(6,3,1)$, and $k=2$. 
\enlargethispage{20mm}

\bigskip
T. Fenske begun the classification of rational cuspidal curves of type
$(d,d-4)$. Recall that a curve $C$ is said to be \textit{unobstructed}
if $H^2(\Theta_V\langle D\rangle)=0$, where 
$(V,D)\rightarrow (\projective, C)$ is a minimal 
embedded resolution of singularities of $C$, 
and $\theta_V\langle D\rangle$ is the sheaf of 
holomorphic vector fields on $V$ tangent along~$D$.
\begin{thm}\label{fe5}\textbf{(Fenske~\cite{Fe2})}
For each $n\geq 1$ there exists a rational cuspidal
plane curve of type $(d,d-4)$, where $d=\mbox{deg}\, C_n=3n+4$. 
This curve has three cusps of types 
$[3n, 3_n]$, $[4_n, 2_2]$, $[2_1]$. The curve 
$C_n$ is rectifiable\footnote{i.e., it is equivalent to a line up to the action of the
  Cremona group ${\bf Bir}(\projective)$ of birational transformations
  of the projective plane.} and unique up
to a projective equivalence. 
Moreover, any unobstructed rational cuspidal
curve of type $(d,d-4)$ is projectively equivalent to a curve of this type.
\end{thm}

\begin{thm}  Groups of the curves in Theorem~\ref{fe5} admit the
  presentation
$$
\fg{\widetilde{C}}=
\langle a, \gamma\,|\, a\gamma a= \gamma a^{n+1} \gamma,\quad
[a^n,\gamma^3]=a^{3n^2+2n}(\gamma^3 a)^{n+1}=1 \rangle.
$$
This group is big provided $3|(n+1)$ and $n>6$. 
\end{thm}
The proof of this theorem is given in 3.3.9.

\section{Calculations}
\textbf{Conventions.}
Throughout this work, we shall use the following conventions:
If $\alpha,\beta:[0,1]\rightarrow T$ are two paths in a topological
space $T$, then the product $\alpha\cdot\beta$ is defined provided 
that $\alpha(1)=\beta(0)$, and one has
$$
\alpha\cdot\beta(t):=
\left\{
\begin{array}{rl}
\alpha(2t), & 0\leq t\leq \frac{1}{2},\\
\beta(2t-1), & \frac{1}{2}\leq t\leq 1.
\end{array}
\right.
$$
If $\alpha$ is a path in $T$ with $\alpha(0)=\alpha(1)=*\in T$, we
shall take the freedom to talk about $\alpha$ as an element 
of $\pi_1(T,*)$, ignoring the fact that the elements of 
$\pi_1(T,*)$ are equivalence classes of such paths under 
the homotopy. Also, when this do not lead to a confusion, 
we shall write $\pi_1(T)$ instead of $\pi_1(T,*)$, omitting the base
point.
  
\subsection{Groups of the curves in Theorem~\ref{fz2}}
\subsubsection{Construction of the curves}

Let $C$ be the cubic defined by the equation $x^2z-y^3=0$.
Then $C$ has a unique, simple cusp
at the point $r=[0:0:1]$, and a unique, simple inflection point at the point   
$p=[1:0:0]$. 
Denote by $P$ the tangent line to $C$ at $p$. 
In order to transform $C$ to $C_n$ by means of
appropriate Cremona transformations,
 we begin by taking an arbitrary point 
$q\in C\moins \{r,p\}$. 
Let $Q$ be the tangent to $C$ at $q$. 
Then $Q$ intersects $C$ at a second point $s$, and the lines $P$ and
$Q$ intersect at a point $O\notin C$ \ffigure{d30}. 

%
%

\figur{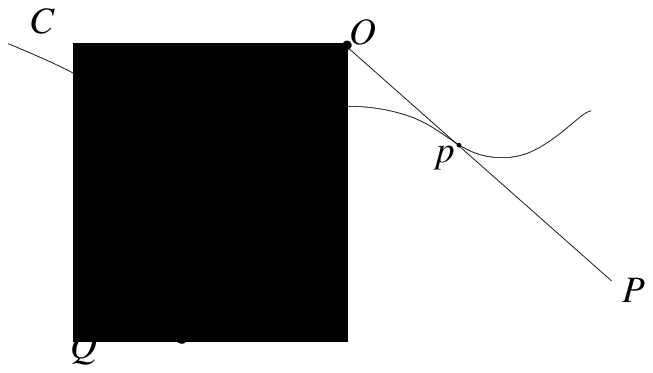}{} 

\par
By blowing-up the point $O$, we obtain a Hirzebruch surface $X$; 
let $E$ be its exceptional section. 
Let $e:=Q\cap E$. 
We apply an elementary transformation (or Nagata transformation) 
at the point 
$e$ followed by an elementary
transformation at the point $s$. Denote by $X_1$ the Hirzebruch
surface so obtained, by $Q_1$, $P_1$ the fibers replacing
$Q$, $P$ and by $\tilde{C}_1$, $E_1$ the proper transforms of $C$, $E$
respectively. 
Then $3p_1:=C_1.P_1$, $2q_1+s_1:=C_1.Q_1$, $e_1:=E_1.Q_1$.
Now we apply an elementary transformation at $s_1$ followed by an 
elementary transformation at $e_1$. 
We get another Hirzebruch surface $X_2$, with fibers $Q_2$, $P_2$ replacing
$Q_1$, $P_1$ and with 
proper transforms $\tilde{C}_2$, $E_2$ of $\tilde{C}_1$, $E_1$
respectively. 
Iterating this procedure $n$ times gives a Hirzebruch surface
$X_n$ with exceptional section $E_n$ satisfying $E_n^2=-1$.
After the contraction of $E_n$, we turn back to $\projective$. 
Then $C_n\subset \projective$ is the  image of $\tilde{C}_n$.

\par
Composition of these birational maps gives a biholomorphism
$$
\projective\moins(C\cup P\cup Q)\stackrel{\simeq}{\longrightarrow}
\projective \moins(C_n\cup P_n\cup Q_n),
$$
and hence an isomorphism
$$
\pi_1(\projective\moins(C\cup P\cup Q))\simeq
\pi_1(\projective\moins(C_n\cup P_n\cup Q_n)).
$$
So, the  group of $C_n$ can be deduced from 
$\pi_1(\projective\moins(C\cup P\cup Q))$
by adding the relations corresponding to 
the gluing of the lines $P_n$ and $Q_n$.

\subsubsection{Finding $\pi_1(\projective\moins(C\cup P\cup Q))$.}
\figur{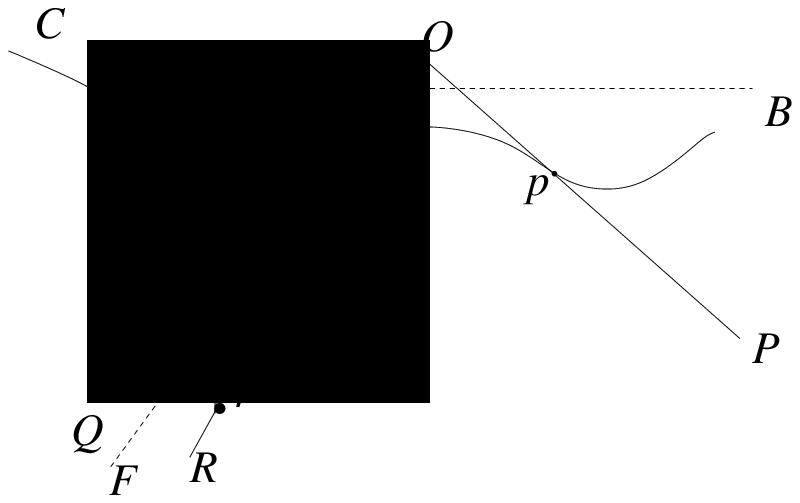}{}

Let $T$ be the projective linear transformation 
$$
[x:y:z]\longrightarrow [x:y:x+z].
$$
Then the equation of $C$ reads, in the new coordinates, as
$x^2(z-x)-y^3=0$, the point $r=[0:0:1]$ is the cusp and $p=[1:0:1]$ is the
inflection point of $C$. 
Put $L_\infty:=\{z=0\}$, and pass to affine coordinates in
$\complex^2=\projective\moins L_\infty$. 
The real picture of $C$ is shown in \fffigure{d34}. 
Let $q=(x_0, y_0)\in C$ be a point such that $x_0>0$ is sufficiently
small, and let $Q$ be the tangent to $C$ at $q$. 
Let $O:=P\cap Q$, and let $R$ be the line $\overline{Or}$.  
Let $B:=\{y=y_1\}$ be a line close to $O$ but $O\notin B$. 
We shall apply the Zariski-Van Kampen method 
to the linear projection $pr:\projective\moins {O}\rightarrow B$ with center $O$. 
Clearly, $P$, $Q$ and $R$ are singular fibers of this projection 
(see \fffigure{d34}).
That these constitute all the singular fibers can be seen by  
looking at the dual picture: 
the dual $C^\star$ of $C$ is known to be the curve $C$ itself (see \cite{Na}). 
The dual $P^\star$ of $P$ is the cusp of $C^\star$, and $O^\star$ is
the tangent line to $C^\star$ at $O^\star$. Since
deg$(C^\star)$=deg$(C)=3$, $O^\star$ cuts $C^\star$ at another point,
which is $Q^\star$. The singular line $R$ corresponds to the
intersection of $O^\star$ with the line $r^\star$.
\par
Consider the restriction $pr^\prime$ of the projection $pr$ to 
$\projective\moins 
(C\cup P\cup Q\cup R)\rightarrow B\moins (P\cup Q\cup R)$. This is a
locally trivial fibration . 
Put $F^\prime:=F\moins (C\cup O)$ where $F$ is a generic fiber of
the projection $pr$, and let $B^\prime:=B\moins(P\cup Q\cup R)$.
Since $\pi_2(B^\prime)=\pi_0(F^\prime)=0$, 
there is the short exact sequence of the fibration 
$$
0\longrightarrow\pi_1(F^\prime)
\longrightarrow\fg{(C\cup P\cup Q\cup R)}
\longrightarrow\pi_1(B^\prime)
\longrightarrow 0.
$$
To determine the group $\fg{(C\cup P\cup Q\cup R)}$ it suffices
therefore to find the \textit{monodromy}, that is, the action of 
$\pi_1(B^\prime)\simeq {\mathbb F}_2$
on the group $\pi_1(F^\prime)\simeq{\mathbb F}_3$.
Choose the base fiber $F$ as shown in \fffigure{d34}. 
Denote by $f_1$, $f_2$, $f_3$ the intersection points $F\cap C$.
Let  $*=F\cap B$ be the base point. 
\par
One can identify the fibers with $\complex$ e.g. by taking $*$ to be
the origin in $\complex^2$, $F$ to be the $y-axis$, and $B$ to be the 
$x-axis$. Then the projection $\complex^2\rightarrow F$ gives the
desired identification. 
\par
Choose positively oriented simple loops  
$a,b,c\in \pi_1(F^\prime,*)$ around $f_1$, $f_2$, $f_3$  
and the loops $\alpha, \beta,\gamma\in \pi_1(B^\prime,*)$ as in
\fffigure{d35}. 
Note that $\pi_1(B^\prime,*)=\langle\alpha, \beta\rangle$.
\figur{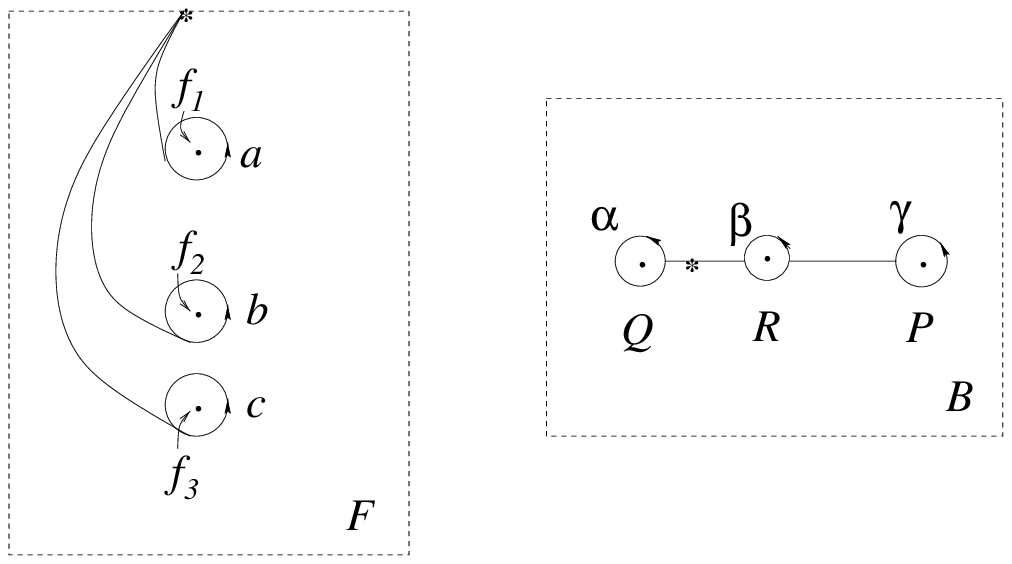}{}
The local monodromy of $pr^\prime$ around the points $p$, $q$, $r$ and
$s$ is well known.
The monodromy around  $q$ gives the relations 
\begin{center}
$\alpha^{-1}a\alpha=b \quad(\rel{1})$\\
$\alpha^{-1}b\alpha=bab^{-1}\quad(\rel{2})$,
\end{center}
and the monodromy around $s$ gives the relation 
$$
\alpha^{-1}c\alpha=c \Leftrightarrow [\alpha,c]=1\quad (\rel{3}).
$$
One has $a=\alpha b\alpha\mon$; by $(\rel{1})$, substituting this in 
the relation $(\rel{2})$ we obtain 
$$
(\alpha b)^2=(b \alpha)^2 \quad (\rel{4}).
$$
The relation obtained from the monodromy around $R$ gives the cusp
relation 
$$
cbc=bcb\quad ({\cal R}_5)
$$
(recall that we glue $R$ back to 
$\projective\moins(C\cup P\cup Q\cup R))$.
Since $\pi_1(B^\prime)=\langle\alpha, \beta\rangle$, 
it is not necessary to calculate the relations obtained 
from the monodromy around $P$; these can be derived 
from the ones we have found.
To sum up, we have the presentation 
$$
\pi_1(\projective\moins(C\cup P\cup Q))=\hspace{9cm}
$$
$$
\langle  a,b,c,\alpha, \gamma\, |\,
a=\alpha b\alpha^{-1},\,
(\alpha b)^2=(b \alpha)^2,\,
cbc=bcb,\,
[\alpha,c]=1,\, 
cba\gamma\alpha=1\rangle,
$$
where the last relation 
$cba\gamma\alpha=1$ 
comes from the loop vanishing at infinity. 
This can be seen as follows: 
Clearly, $cba$ is a loop in $F$
surrounding the points $f_1$, $f_2$, $f_3$. 
Let $\Sigma$ be a small disc in $F$ containing the point $O$, and let
$\sigma$ be its boundary. 
Let $\realline$ be the real line in $F$, put $h:=\sigma\cap
\realline$, and let $\omega$ be the real line segment $\overline{*h}$. 
Define the positively oriented loop $\rho$ as
$$
\rho:=\omega\cdot\sigma\cdot\omega\mon.
$$
Then one has the relation $\rho cba=1$ since $F=\sphere=S^2$ is a sphere.
Let $U$ be a small neighborhood of $O$
in $\projective\moins (C\cup P\cup Q)$. 
Then clearly $U$ is biholomorphic to $\Delta^*\times \Delta^*$, where 
$\Delta^*$ is the punctured disc (see \fffigure{aroundO}).
Hence, 
$\pi_1(U)={\integers}^2=\langle\alpha, \gamma\,|\, [\alpha,\gamma]=1\rangle$,  
and it is easy to see that $\rho$ is homotopic to $\alpha\gamma$. 
\figur{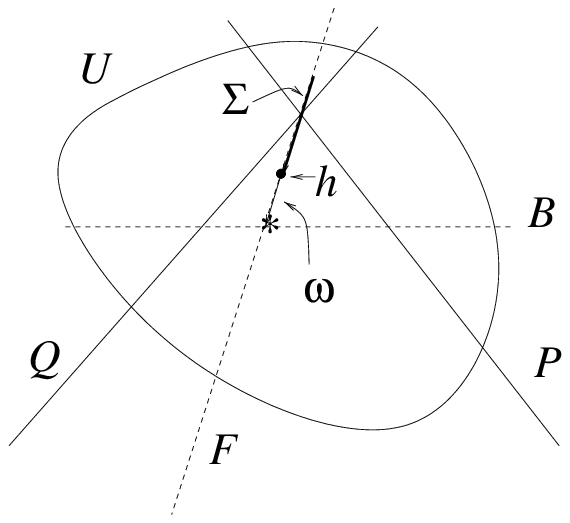}{}
\par
Using $({\cal R}_1)$ and $[\alpha,\gamma]=1$, the relation
$cba\gamma\alpha=1$ becomes
$$
cb\alpha b\gamma=1\quad ({\cal R}_6).
$$
Eliminating the generators $a$ and $\gamma$
from this presentation, we get 
$$
\pi_1(\projective\moins (C\cup P\cup Q))=
\langle b,c,\alpha\,|\,
(\alpha b)^2=(b \alpha)^2,\,
cbc=bcb,\,
[\alpha,c]=1\rangle
$$
$$
=\pi_1(\projective\moins(C_n\cup P_n\cup Q_n)).
$$
To obtain a presentation of the group $=\pi_1(\projective\moins C_n)$,
it remains to find the relations corresponding to the 
gluing of the lines $P_n$ and $Q_n$. To this end, we introduce the
following concept.

\begin{defin}
\label{definmeridian}
\textbf{(meridian)}
\rm 
Let $C$ be a curve in a surface $X$, 
and pick a base point $*\in X\moins C$. 
Let $\Delta$ be a small analytic disc in $X$, intersecting $C$
transversally at a unique  point $q$ of $C$.
If $q$ is a smooth point of $C$, a
 \textit{meridian} of $C$ in $X$ with respect to the base point $*$
is a loop in $X\moins C$ constructed 
as follows:
\figur{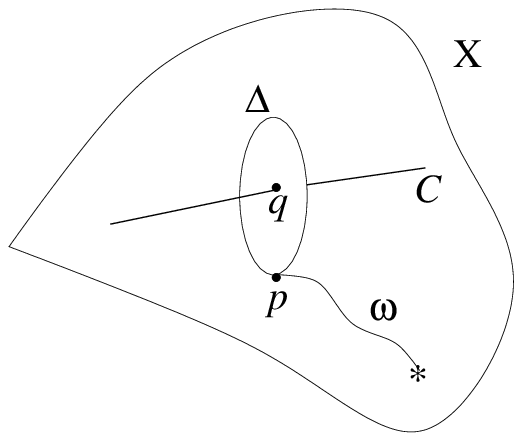}{}
 Connect $*$ to a point $p\in\partial\Delta$ by means
of a path $\omega\subset X\moins C$ such that 
$\omega\cap\Delta=p$, and let
$$
\mu:=\omega\inv \cdot \delta\cdot \omega,
$$  
where $\delta:=\partial\Delta$, oriented clockwise
(\fffigure{meridien}). A loop $\mu_q$ given by the same construction will be
called a \textit{singular meridian at $q$} if $q$ is a singular point of $C$.
\end{defin}
\begin{lem}
\label{meridians}
Let $M\vartriangleleft  \pi_1(X\moins C,*)$ be the 
subgroup normally generated by the meridians
of $C$. Then 
\\(i) $\pi_1(X)=\pi_1(X\moins C,*)/M$. 
\\(ii) If $C$ is irreducible, then any two meridians of $C$ are
conjugate in $\pi_1(X\moins C,*)$.
 Hence, $M=\,\ll\mu \gg$, 
where\footnote{Recall that by $\ll\mu \gg$ we 
 denote the normal subgroup generated by $\mu$.}
 $ \mu $ is any meridian of $C$.
\par(iii) Any two singular meridians $\mu_q$, $\tilde{\mu}_q$ at the
singular point $q\in C$ are conjugate.
\end{lem}   
\proof Parts \textit{(i)-(ii)} are well known \cite{Lamotke}. 
To show \textit{(iii)}, 
assume that $\mu_q$, $\tilde{\mu}_q$ are
obtained from the discs
$\Delta_q$, $\widetilde{\Delta}_q$ intersecting $C$
transversally at $q$.  Let $\sigma:\,Y\rightarrow X$ be the blow-up of the
surface $X$ at the point $q$, and denote by $Q:=\sigma\inv(q)$ 
the exceptional divisor of this blow-up. Then the proper transforms 
$\sigma\inv(\Delta_q)$, $\sigma\inv(\tilde{\Delta}_q)$ 
intersect $Q$ transversally at
distinct points of $Q$, and these points of intersection are smooth in
$C\cup Q$. Hence, $\sigma\inv(\mu_q)$, $\sigma\inv(\tilde{\mu}_q)$ are
meridians of $Q$. 
Since $Q$ is irreducible, applying the part \textit{(ii)} to the
surface $Y\moins C$ gives the desired result. $\Box$

\bigskip

\figur{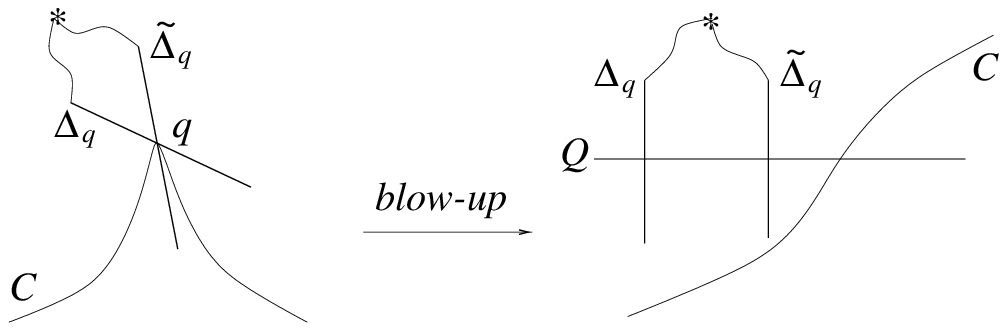}{}
\par
For a group $G$, denote by $(g)$ the conjugacy class of
$g\in G$ , i.e. $(g):=\{hgh\mon\,:\, h\in G\}$. 
Lemma~\ref{meridians} implies that the  group of a curve
$C\subset \projective$,
supplied with the following data
\\\textit{(i)} Conjugacy classes $(\mu_1), \,(\mu_2),\dots$ 
of meridians of $C$,
\\\textit{(ii)} Conjugacy classes of singular meridians 
$(\mu_{q_1}), \,(\mu_{q_2}),\,\dots$ of $C$ at singular points
$q_1,\,q_2,\dots$ of $C$ \\
is a richer invariant of the pair  $(\projective, C)$ 
than solely the group $\fg{C}$. 
\par
For the group $\pi_1({\projective\moins (C\cup P\cup Q}))$ found above, 
there are clearly three classes of meridians $(a), \,(\alpha),
\,(\beta)$, corresponding to the curves $C,\,P,\,Q,$ respectively.
Consider the loop $cb$ in $F$, which surrounds the points 
$f_2$ and $f_3$. Pushing $F$ over $R$, the points $f_2$, $f_3$ come
together at the cusp $r$, and the loop $cb$ becomes a loop in $R$ 
surrounding the cusp $r$, that is, $cb$
is a singular meridian of $C\cup P\cup Q$ at $r$, i.e. $(\mu_r)=(cb)$. 
(The classes $(\mu_O)$, $(\mu_p)$, $(\mu_q)$, $(\mu_s)$ are irrelevant
so we will not find them.)
\par
On the other hand, Lemma~\ref{meridians}
 implies that the relations corresponding to 
the gluing of the lines $Q_n$,
$P_n$ are of the form $\mu(Q_n)=1$, $\mu(P_n)=1$,
where $\mu(Q_n)$, $\mu(P_n)$ are meridians of 
$Q_n$ and $P_n$ respectively. 
Finding these meridians will be achieved
by  Fujita's lemma, which we proceed to explain now.
\\
Let  $q\in C$ be an ordinary double point, and take
a small neighborhood $X'$ of $q$ such that $X'\cap C$ consists
of two branches $C_1$ and $C_2$ satisfying $C_1\cap C_2=\{q\}$.  
Pick an intermediate base point $*'\in X'\moins C$, and take meridians
$\mu_1'$ of $C_1$ in $X'\moins C_2$ and $\mu_2'$ of $C_2$ in 
$X'\moins C_1$ with respect to $*'$. 
Let $\omega$ be a path in $X\moins C$ connecting $*$ to $*'$, and
define
$$  
\mu_1:=\omega\inv \cdot \mu_1'\cdot \omega,\quad
\mu_2:=\omega\inv \cdot \mu_2'\cdot \omega.
$$ 
Clearly, $\mu_1$, $\mu_2$ are meridians of $C$ in $X$ with respect to
$*$ and they commute. 
Moreover, $\mu_1\cdot\mu_2$ is homotopic to a 
singular meridian of $C$ at $q$.
Hence, we have the following lemma:
\figur{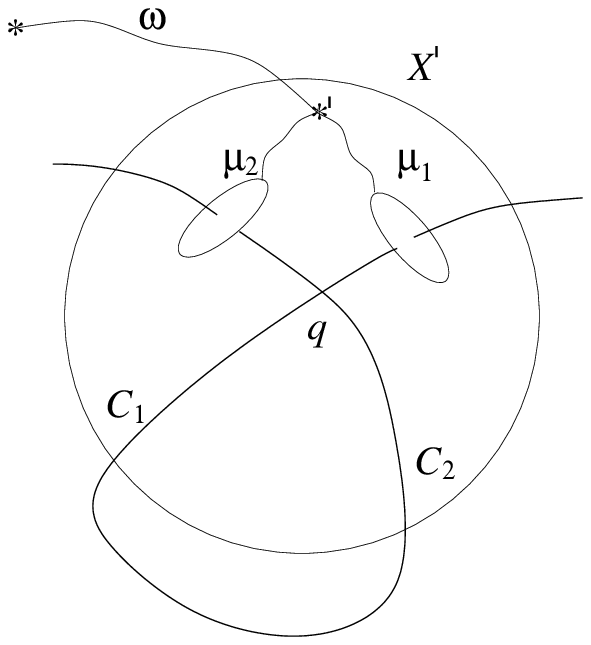}{}
\begin{lem}\textbf{(Fujita~\cite{Fu})}\label{fujita1}
Let $\sigma: X\mapsto Y$ be the blowing up of $q$, and put
$Q:=\sigma^{-1}(q)$. Identify $X\moins \{Q\}$ with $Y\moins \{q\}$.
Then $\mu_1\cdot\mu_2$ is a meridian of $Q$ in $Y$ with respect to
$*$. Moreover, as $Q$ is irreducible, by Lemma~\ref{meridians} one has 
$$
\pi_1(X\moins C,*)=
\pi_1\left(Y\moins (C\cup Q),*\right)/\ll\mu_1\cdot\mu_2\gg.
$$
\end{lem} 
\subsubsection{Meridians of $P_n$ and $Q_n$}
Turning back to our search for the relations introduced in 
$\pi_1(\projective\moins(C_n\cup P_n\cup Q_n))$ after the gluing of
$P_n$, $Q_n$, we first note that one can apply Fujita's Lemma to the
loops $\alpha$, $\gamma$, which are meridians of $P$ and $Q$
respectively, and obtain a meridian of $E$. The blowing-up of the
point $E\cap P$ will give $\gamma(\alpha\gamma)$ as a meridian of $P_1$,
by induction we obtain 
$\mu(P_n):=\gamma(\alpha\gamma)^n$ as a meridian of $P_n$. 
Recalling that $[\alpha,\gamma]=1$, this gives the relation
$$
\alpha^n\gamma^{n+1}=1.
$$ 
Substituting $\gamma$ from $({\cal R}_6)$ we get,
$$
\alpha^n (cb\alpha b)^{-(n+1)}=1\quad ({\cal R}_7).
$$
\par
The construction of $C_n$ also shows that 
$\mu_p:=\alpha^{n-1}\gamma^n$ is a singular meridian of $C_n$ at
$p_n$, since it is a meridian of the exceptional line of the blow-up
at $p$. Setting $\mu(P_n)=1$ yields that $\mu_p=(\alpha\gamma)\mon$
is a singular meridian of $C_n$ at $p_n$. 
\figur{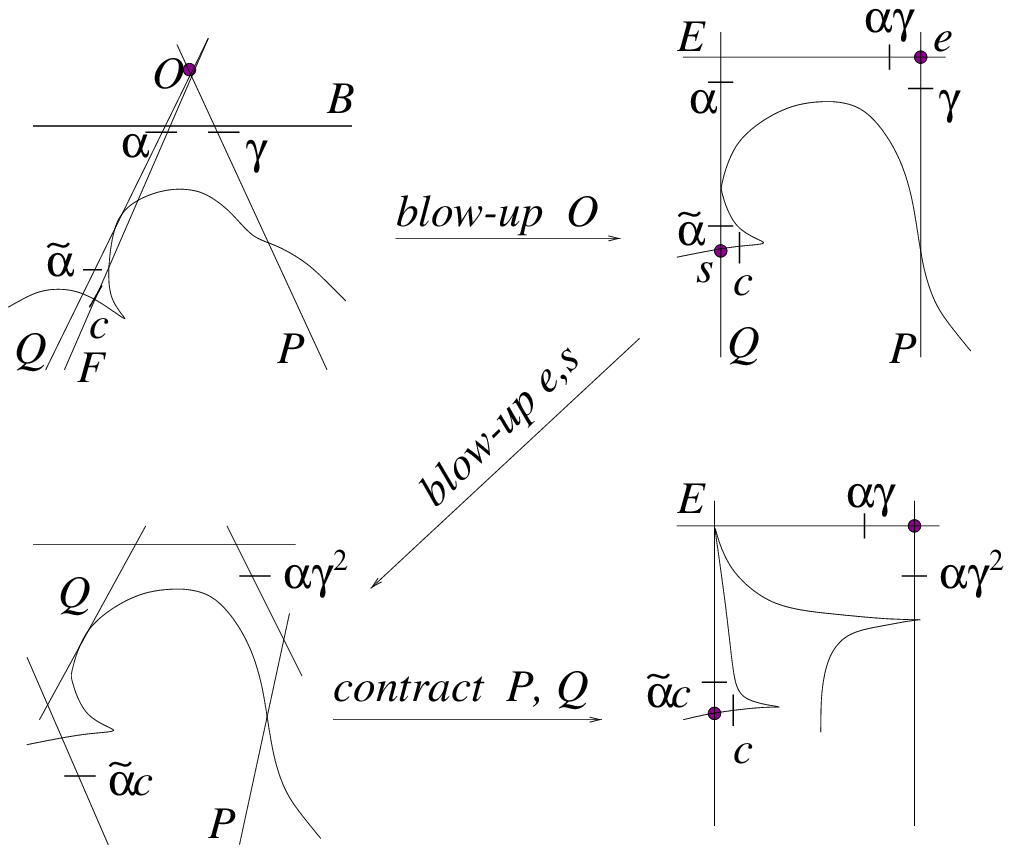}{}
\par
To find a meridian of $Q_n$, 
first define  a meridian $\tilde{\alpha}$ of $Q$ as
shown in figure~(\ref{meridiens}). 
That is, take a small disc $\Delta$
intersecting $Q$ transversally above the point $s$, and take a path 
$\omega$ joining $\Delta$ to a neighborhood $f_3$, 
and continuing to $*$ in $F^\prime$ along the loop $c$. 
Let $\tilde{\alpha}:=\omega\cdot\delta\cdot\omega\inv$. 
Then the blowing up of the point 
$s$ will give $\tilde{\alpha}c$ as a meridian of $Q_1$. 
A recursive application of Fujita's Lemma gives 
$\mu(Q_n):=\tilde{\alpha}c^n$ as a meridian of $Q_n$. 
\par
Note that the construction of $C_n$ also shows that 
$\mu_q:=\alpha\gamma$ is
a singular meridian of $C_n$ at $q_n$, since $\mu_q$
is a meridian of the exceptional line of the blow-up of $q_n$.
\par
\subsubsection{Finding $\tilde{\alpha}$}
The final step in determining the fundamental group of $C_n$ is to
express $\tilde{\alpha}$ in terms of the above presentation of 
$\pi_1(\projective\moins(C\cup P\cup Q))$. Let us show that
$\tilde{\alpha}$ is in fact homotopic to $\alpha$. 
\par
First, let $\widetilde{\Delta}:=pr\mon(\Delta)\cap B$. 
Then $\alpha$ is clearly homotopic to a loop obtained by connecting 
(properly) $\partial \widetilde{\Delta}$ to the base point $*$.
$\partial \widetilde{\Delta}$ intersect the real axis of $B$ at two
points, let $*^\prime$ be the one on the right, which will be used as a
temporary base point. Now push the fiber $F$ over $*^\prime$ along the
real axis of $B$. As all the intersection points remains real, it is
easy to see that the picture of $F$ stays as in \fffigure{d35}. 
\par
Next, consider the restriction 
 $pr': {pr^{-1}(\partial\widetilde{\Delta})}\rightarrow 
\partial\widetilde{\Delta}$ 
of $pr$ to the border of the disc $\widetilde{\Delta}$. 
This is a locally trivial fibration, and it can be pictured as in  
\fffigure{d3dragalfa3}, where we have cut 
$\partial\widetilde{\Delta}$ at $*^\prime$ to give a better picture. 
Let $\Sigma$ be a disc in $F$,  
containing the upper intersection
points $f_1$ and $f_2$ of $C$ with $F$, we suppose that $\Sigma$
avoids the loop $c$. Then, as the lower intersection point $f_3$ is a
transversal intersection, we can suppose that the corresponding 
Leftschez homeomorphisms $F\rightarrow F_t$, 
$t\in\partial\widetilde{\Delta}$ are constant outside $\Sigma$. 
Hence, the following map gives a homotopy between $\alpha$ and 
$\tilde{\alpha}$.
\figur{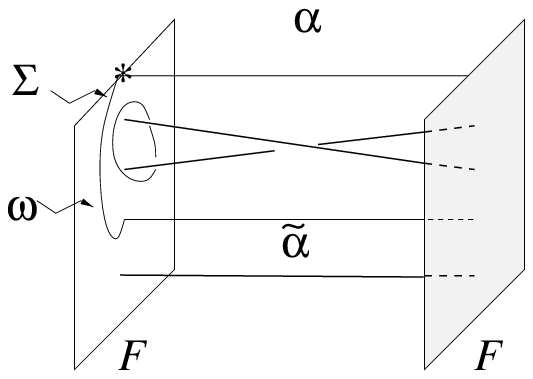}{}
$$
H(s,t):=
\left\{
\begin{array}{rl}
(\alpha(0), \omega(t)), & 0\leq t\leq s/3\\
(\alpha\left(\frac{t-s/3}{1-2s/3}\right), \omega(s/3)),& s/3\leq t\leq 1-s/3\\
(\alpha(1), \omega(1-t)) & 1-s/3\leq t\leq 1
\end{array}
\right.
$$
Consequently, the vanishing of the meridian of $Q_n$ yields the relation
$$
\mu(Q_n)=1\Rightarrow\alpha=c^{-n}\quad({\cal R}_8).
$$
Substituting $\alpha$ from $({\cal R}_8)$, $({\cal R}_7)$ becomes
$$
c^{n^2}(cbc^{-n}b)^{n+1}=1.
$$
From the cusp relation $cbc=bcb$ it follows that 
$cbc^{k}=b^{k}cb$ for any $k\in\integers$. 
Using this in the above relation, we get
$$
c^{n^2}(b^{-n}cb^2)^{n+1}=1 \quad ({\cal R}_9).
$$
To sum up, we have the presentation
$$
=\pi_1(\projective\moins C_n)
=\langle b,c,\alpha \,|\, 
{\cal R}_3,\,
{\cal R}_4,\,
{\cal R}_5,\,
{\cal R}_7,\,
{\cal R}_9 \rangle.
$$
\subsubsection{Study of the group}
Since $\alpha=c^{-n}$ by $({\cal R}_7)$, the relation $({\cal R}_3)$ is
trivialized, and using the cusp relation $({\cal R}_4)$ becomes
$$
(bc^{-n})^2=(c^{-n}b)^2\Leftrightarrow [b^n,c^{n+2}]=1.
$$
Finally, we have obtained the presentation given in 
Theorem~\ref{groupssflza}.
$$
G_n:=\pi_1(\projective\moins C_n)=
\langle b,c \,|\, cbc=bcb,\quad [b^n, c^{n+2}]=1,\quad
(b^{-n}cb^2)^{n+1}c^{n^2}=1\rangle.
$$
\par
Notice that $(c)$ is the unique class of meridian in $G_n$. Also,
the singular meridian $\mu_r=cb$ of $C$ is unchanged during the birational
transformations, so that $\mu_r$ is a singular meridian of $C_n$ at $r_n$.
Other singular meridians are, as found above, $\mu_q=\alpha\gamma$ and 
$\mu_p=(\alpha\gamma)\mon$. 
One has $(\alpha\gamma)\mon=cba=cb\alpha b\alpha\mon=cbc^{-n}bc^n$.
\\
It is easily seen that $G_0=\integers/3\integers$. 
Let us now show that $G_1=\integers/5\integers$ and 
$G_2=\integers/7\integers$. One has
$$
G_1=\langle c,b,\,|\, cbc=bcb, \, [b,c^3]=1,\, 
(b\mon cb^2)^2c=1\rangle.
$$
Expanding the last relation, we get
$$
b\mon c\cdot bcb^2\cdot c=
b\mon c\cdot c^2bc\cdot c=
b\mon c^3bc^2=c^5.
$$
Thus,
$$
[b,c^3]=1\Rightarrow [b,c^6]=[b,c\mon]=[b,c]=1\Rightarrow b=c.
$$
Hence, $G_1$ is generated by $b$, and $G_1=\integers/5\integers$. The
group $G_2$ has the presentation
$$
G_2=\langle c,b\,|\, cbc=bcb,\quad 
[b^2,c^4]=1,\quad (b\mtw cb^2)^3c^4=1\rangle.
$$
Again, expanding the last relation we get
$$
b\mtw c^3 b^2 c^4=1\Rightarrow c^7=1.
$$
Thus,
$$
[b^2, c^4]=1\Rightarrow [b^6, c^8]=[b\mon, c]=[b,c]=1\Rightarrow c=b.
$$
Hence, $G_2$ is generated by $b$, and one has
$G_2=\integers/7\integers$.
The group $G_3$ is found to be finite of order $8640$ 
by using the programme Maple. The order of the group $G_5$
is calculated to be $1560$ by Artal by the help of the programme GAP.

\begin{defin}\label{definresidual}\textbf{(residual group)}
\rm
 For an element $a\in G$ of an arbitrary
group $G$, the group $G/\ll a\gg$ will be denoted by $G(a)$. 
If $G_C:=\fg{C}$ is the  group of an
irreducible curve $C$, with $\mu$ a meridian of $C$, 
then for $k\in \naturals$, we will call the group $G_C(\mu^k)$ 
a \textit{residual group} of $G_C$ and denote it by $G_C(k)$. 
The group $G_C(\mu_p^k)$, where $\mu_p$ is a singular meridian of
$C$ at a singular point $p\in C$, will be called a 
\textit{residual group of} $G_C$ \textit{at} $p$ and denoted by $G_C^p(k)$. 
Note that the groups $G_C(k)$  
do not depend on the particular meridian $\mu$  
chosen, since by Lemma~\ref{meridians}(ii),
for an irreducible curve, any two meridians  
are conjugate. In view of Lemma~\ref{meridians}(iii),
this is also true for the groups $G_C^p(k)$.
\end{defin}
\begin{prop}
\label{trianglegroup}
If $n$ is odd, and $k|n$, then there is a surjection
$$
G_n(k)\onto T_{2,3,k}=\langle x,y,z\,|\, x^2=y^3=z^k=xyz=1\rangle
$$
onto the triangle group $T_{2,3,k}$. Hence, $G_n$ is big for 
odd $n\geq 7$.
\end{prop}
\proof
For the last assertion,
it is known that the group $T_{2,3,k}$ 
is big if $k\geq 7$.
Putting $n=k$, we get that the groups 
$G_n(n)$ are big for $n\geq 7$ odd, so that $G_n$ is big
if $n\geq 7$ is odd.
\par
Now let us establish the surjection claimed. 
If $b^k=1$, then $c^n=b^n$=1 since $k|n$. This gives the presentation 
$$
G_n(k)=
\langle b,c \,|\, 
cbc=bcb,\quad
(cb^2)^{n+1}=1,\quad
b^k=1\rangle.
$$
Following an idea due to Artal \cite{Artal1}, we 
apply the  transformation $x=cbc$, $y=cb$
($\Leftrightarrow$ $c=y\mon x$, $b=x\mon y^2$) to obtain
$$
G_n(k)=
\langle x,y \,|\,
x^2=y^3,\quad
 (yx\mon y^2)^{n+1}=1,\quad
(y\mon x)^k=1\rangle.
$$
Note that $yx\mon y^2=yxy\mon$ since $x^2=y^3$. Thus,
$$
G_n(k)=
\langle x,y \,|\,
x^2=y^3,\quad
x^{n+1}=1,\quad
(y\mon x)^k=1\rangle.
$$
Let $H$ be the quotient of this group
by the relation $x^2=y^3=1$ (note that $x^2$ is central).
Then the relation $x^{n+1}=1$ is
killed if $n$ is odd, and we get the desired result: 
$$
H=
\langle x,y \,|\,
x^2=y^3=(y\mon x)^k=1\rangle= T_{2,3,k}.\qquad \Box
$$
This completes the proof of Theorem~\ref{groupssflza}.

\subsection{Groups of the curves (1)-(1a) in Theorem~\ref{table}}
\subsubsection{Construction of the curves}
Before passing to the construction of the curves, let us fix some
notations following \cite{Fe1}.

\bigskip\noindent\textbf{Notation.}
 Let $\sigma_O: X_1\rightarrow \projective$ be the blow-up of 
$\projective$ at the point $O\in\projective$. 
Denote by $E_1$ the exceptional divisor of this blow-up.  
For any curve $C\subset \projective$, the proper preimage of $C$
in the Hirzebruch surface $X_1$ will be denoted by $C_1$. 
Now let $X_i$ be a Hirzebruch surface. 
Then $X_i$ is a ruled surface, whose horizontal section is denoted by
$E_i$. Let $p\in X_i$, and let $L_i$ be the fiber of the ruling
passing through $p$. 
The surface obtained from $X_i$ by an elementary transformation at the
point $p$ will be denoted 
by $X_{i+1}$. Recall that this is a birational mapping which consists
of blowing-up  $p\in X_i$ followed by the 
contraction of $L_i$. The fiber replacing $L_i$ will be denoted by 
$L_{i+1}$, and for any curve $C_i\in X_i$, the proper transform of 
$C_i$ in $X_{i+1}$ will be denoted by $C_{i+1}$. 

\bigskip
Let $n,m\in \naturals$ be such that $0\leq m< n$, 
and for $k\geq 2$, let
$C$ be the curve defined by the equation $F(x,y,z):=xy^{k-1}-z^k=0$.  
Then for $k>2$, the curve 
$C$ has a unique singularity at $p:=[1:0:0]$ which is a cusp and 
$q=[0:1:0]$ is an inflection point of $C$ of order $k$. 
If $k=2$, then $p$ is a smooth point of $C$.
The line $P:=\{y=0\}$ is the tangent to $C$ at $p$ and  
$Q:=\{x=0\}$ is the tangent at $q$, these tangents 
intersects at the point $O:=P\cap Q=[0:0:1]\notin C$. 
Blowing-up the point $O$, we get a Hirzebruch surface $X_1$. 
Let $E_1$ be its horizontal section, and denote by $P_1$, $Q_1$
the proper transforms of $P$ and $Q$ respectively. 
For $i=1, 2, \dots, m$ we apply $m$ 
elementary transformations at the points  
$E_i\cap P_i$, followed by elementary 
transformations applied at the points 
$E_i\cap Q_i$ for $i=m+1,m+2,\dots, n$,
and we arrive at the Hirzebruch surface $X_{n+1}$ with 
$E_{n+1}^2=-n-1$ (see \fffigure{fenske3}).
\par
Performing elementary transformations at arbitrary points 
$s_i\in P_i\moins E_i$ for $i=n+1, \dots, 2n$ we obtain the Hirzebruch 
surface $X_{2n+1}$ with $E_{2n+1}^2=-1$. 
Hence, one can contract $E_{2n+1}$ and return to the projective 
plane $\projective$. 
Let $\widetilde{C}, \widetilde{P}, \widetilde{Q}$ be the images
 of respectively $C_{2n+1}, P_{2n+1}, Q_{2n+1}$ 
under the contraction of $E_{2n+1}$. 
Then  $\widetilde{C}$ is a curve of the family (1).
\par
The curves (1a) are obtained in the same way, except that in this case
 one applies elementary transformations at the points 
 $E_i\cap P_i$ for  $i=1,2,\dots,n$ followed by 
elementary transformations at some points 
 $s_i\in E_i\cap Q_i$ for  $i=n+1,n+2,\dots,2n$.
\par
These birational morphisms provides a biholomorphism
$$
\projective\moins(C\cup P\cup Q)
\stackrel{\simeq}{\longrightarrow}
\projective\moins(\widetilde{C}\cup \widetilde{P}\cup \widetilde{Q}).
$$ 
One has the induced isomorphism
$$
\pi_1(\projective\moins(C\cup P\cup Q))
\simeq
\pi_1(\projective\moins(\widetilde{C}\cup \widetilde{P}\cup \widetilde{Q})).
$$ 
\subsubsection{Finding $\pi_1(\projective\moins(C\cup P\cup Q))$} 
We will apply the Zariski-Van Kampen method to the projection 
$pr:=\projective\moins {O}\rightarrow \sphere$. 
Clearly, $P$, $Q$ are singular fibers of this projection, and it is
easy to see that these are the only ones. 
Indeed, a line passing through $O=[0:0:1]$ has an equation of the form 
$ax+by=0$. Comparing with the equation $xy^{k-1}-z^k=0$ of $C$, 
one obtains $by^k+az^k=0$, which has multiple solutions if and only if 
$a=0, b\neq 0$ or $a\neq 0, b=0$, corresponding to the lines 
$P=\{y=0\}$ and $Q=\{x=0\}$.  Let $pr^\prime$ be the restriction of
the projection $pr$ to $\projective\moins (C\cup P\cup Q)$.
\par
\figur{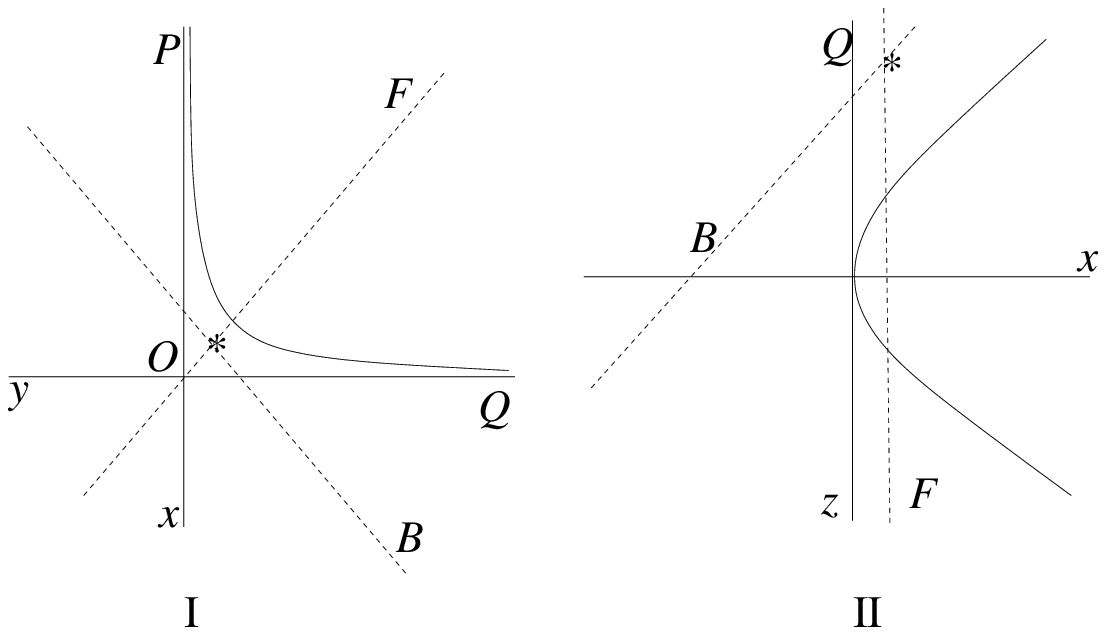}{}
Let $L_\infty:=\{z=0\}$, and shift to the affine coordinates in
$\projective\moins L_\infty =\complex^2$.
Let $B$ be the line $\{x+y=\epsilon\}$, where $\epsilon$ is a small
real number and let $F:=\{x=y\}$ be the base fiber (\fffigure{C_k1}-I). 
Put $B^\prime:=B\moins(P\cup Q)$ and $F^\prime:=F\moins (\{O\}\cup C)$. 
If we choose $L_\infty:=\{y=0\}$ and pass to the affine coordinates in
$\projective\moins L_\infty =\complex^2$, then the real picture 
of the configuration $C\cup P\cup Q$ will be as it is 
drawn in \fffigure{C_k1}-II.
Let $*:=F\cap B$ be the base point. Identify the fibers with
$\complex$ via the projection to the $z$-axis, 
and take the generators of $\pi_1(B^\prime)$ and $\pi_1(F^\prime)$ 
as in \fffigure{C_k}.
\figur{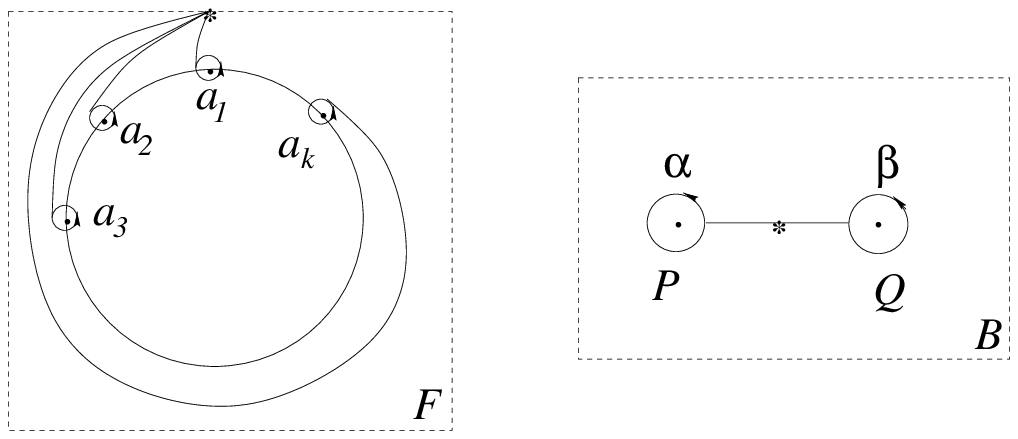}{}
The monodromy relations around the singular fiber $Q$ are given by 
$$
\beta^{-1} a_i \beta =
\left\{
\begin{array}{ll}
a_{i+1},& 1\leq i\leq k-1,\\
      \delta \;a_1\; \delta^{-1}, & i=k.
\end{array}
\right.
$$
where $\delta:=a_k a_{k-1}... a_2 a_1$.
Setting $a:=a_1$, these relations can be expressed as 
$$
\label{fensketangentrelations}
a_i=\beta^{-i+1}a\beta^{i-1}, \quad (\beta a)^k=(a\beta)^k.
$$
Hence, one has the presentation
$$
\pi_1(\projective\moins(C\cup T\cup L))=
\langle\beta, a \quad|\quad (a\beta)^k=(\beta a)^k\rangle.
$$
Note that $[\alpha,\beta]=1$ and $\alpha\beta a_k\cdots a_2a_1=1$.
The change of generators 
$(\beta,\quad a)\Leftrightarrow (\beta,\quad y:=\beta a)$ 
gives a more convenient presentation 
$$
\pi_1(\projective\moins(C\cup T\cup L))
=\langle\beta,y \;|\; [\beta, y^n]=1\rangle.
$$
For the future applications, note that $\alpha$ can be expressed
by using the relation
$$
\alpha\beta a_ka_{k-1}\cdots a_1=1\Leftrightarrow \alpha=y^{-k}\beta^{k-1}.
$$
Note also that $[\alpha,\beta]=1$. This can be derived either from the
above presentation or by applying Fujita's lemma to the meridians 
$\alpha$ of $P$ and  $\beta$ of $Q$, with respect to the point 
$O=P\cup Q$. 
\subsubsection{Meridians of $\widetilde{P}$ and $\widetilde{Q}$}
\figur{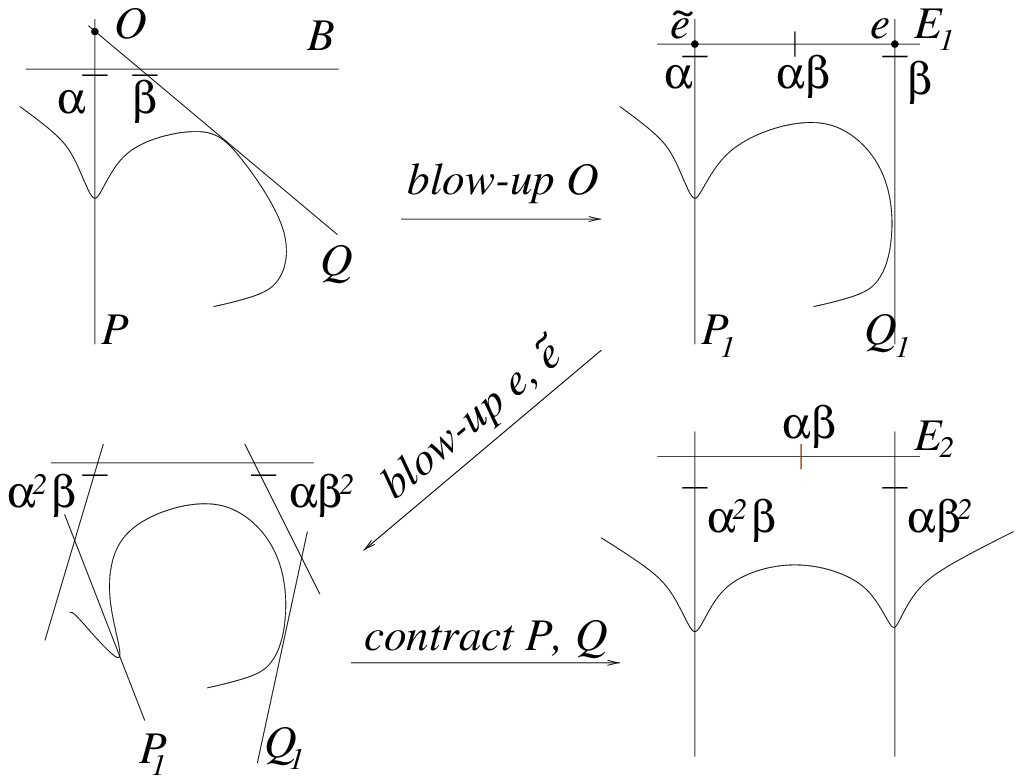}{}
An obvious application of Fujita's lemma yields that 
$\mu(P_{n+1}):=\alpha(\alpha\beta)^m$ is a meridian of $P_{n+1}$ and
$\mu(Q_{n+1}):=\beta(\alpha\beta)^{n-m}$ is a meridian of $Q_{n+1}$ 
in the surface
$X_{n+1}\moins (C_{n+1}\cup P_{n+1}\cup Q_{n+1}\cup E_{n+1})$.
(See \fffigure{fenske3}, where the situation is illustrated for $n=2$, 
$m=1$, beware that the elementary transformation applied at 
$e:=E_1\cap P_1$ and the elementary transformation applied at 
$\tilde{e}:=E_2\cap Q_2$ are shown simultaneously in the figure.)
Recall that the subsequent transformations are applied at points 
$s_i\in E_i\moins P_i$ for $i=n+1,\dots ,2n$. 
Since the line $Q_{n+1}$ is not affected by these
transformations, $\mu(Q_{n+1})$ is a meridian of $\widetilde{Q}$ in  
$\projective\moins(\widetilde{C}\cup 
\widetilde{P}\cup \widetilde{Q})$, too. 
\par
On the other hand, $\mu(P_{n+1})$ stays to be a 
meridian of $P_{n+i}$ after an elementary 
transformation applied at a point $s_i\in E_i\moins P_i$. 
This can be seen e.g. by choosing $s_{n+1}=\Delta\cap P_{n+1}$, 
where $\Delta$ is the defining disc of $\mu(P_{n+1})$.
Hence, $\mu(P_{n+1})$ is a meridian of $\widetilde{P}$, too.
\par
Denote by $\tilde{p}$ and $\tilde{q}$ the cusps of 
$\widetilde{C}$. Then, by the construction of the curve, 
$\alpha\beta$ is a singular meridian at $\tilde{p}$, and 
$\beta(\alpha\beta)^{n-m-1}$ is a singular meridian at 
$\tilde{q}$.
\par
Setting $\mu(P_{n+1})=\mu(Q_{n+1})=1$, we obtain the presentation
$$
G_{(1)}:=\fg{\widetilde{C}}=
$$
$$
\langle \alpha, \beta, y\,|\,
\alpha=y^{-k}\beta^{k-1},\quad
[\beta,y^k]=\alpha(\alpha\beta)^m=\beta(\alpha\beta)^{n-m}=1
\rangle.
$$
\par
A meridian of $\widetilde{C}$ can be given as $a=\beta\mon y$.
After the obvious simplifications, one finds that 
$\alpha\beta=y^{-k}\beta^k$ is a singular meridian at $\tilde{p}$,
and $(\alpha\beta)\mon= y^k\beta^{-k}$ is a singular meridian at 
$\tilde{q}$.
\par
It is easy to see that the element $y^k$ is central in the group $G_{(1)}$. 
Eliminating $\alpha$ in $G_{(1)}(y^k)$, one obtains
$$
G_{(1)}(y^k)=\langle \beta, y\,|\,
y^k=\beta^j=1 \rangle=\integers_k*\integers_j,
$$
where $j:=\mbox{g.c.d.}(mk+k-1, nk-mk+1)=
\mbox{g.c.d.}(mk+k-1, kn+k)=\mbox{g.c.d.}(mk+k-1, n+1)$. 
Hence, this group is big if $j \geq 2$. 
Finally, the group $G_{(1)}$ is 
abelian when $j=1$ by the following trivial lemma:
\begin{lem}
Let $G$ be a group, and $z\in G$ be a central element.
If $G(z)$ is cyclic, then $G$ is abelian.
\end{lem}
\par
As for the curves from the family (1a), the same procedure applies. 
The meridian $\beta$ of $Q$ stays to be a meridian of
$\widetilde{Q}$, so that one has the relation $\beta=1$, 
which implies that the fundamental group is generated by 
just one element $y(=a)$ and thus it is abelian. 

\bigskip\noindent\textbf{Remark.}
In~\cite{DOZ}, the authors provide  a long argument due to
V. Lin, showing the bigness of the group given by the presentation
$\langle a,b\,|\, (ab)^2=(ba)^2\rangle$. Here is a simpler proof of a
more general assertion:
\begin{prop}\label{big}
Let $n,m\in \integers$ such that $k:=gcd\,(n,m)$ satisfy $|k|\geq 2$. 
Then the group
$$
\langle a,b\,|\, (ab)^n=(ba)^m\rangle
$$
is big.
\end{prop}
\proof
Put, as above, $x:=ab$ and $y:=b$. Then the above presentation is
written in terms of $x,y$ as
$$
\langle x,y\,|\, x^n=yx^my\mon\rangle.
$$
Passing to the quotient by the relation $x^k=1$, we get the group 
$\integers_k*\integers$, which is big. This can be seen as follows:
Let $r\geq 3$ be an integer such that $gcd\,(r,k)=1$. 
Passing once more to the quotient by the relation $y^r=1$ gives the
group $\integers_k*\integers_r$, and it is well known that the
commutator subgroup of this group is the free group of rank 
$(k-1)(r-1)$. $\Box$

\bigskip\noindent
Note that the group $\integers_2*\integers_2$ is isomorphic to the
infinite dihedral group ${\mathbb D}_\infty$, whose commutator 
subgroup is  $\integers$, hence this group is solvable and is  not
big.  
Also, it can be shown that the commutator subgroup of the group 
$\langle a,b\,|ab=(ba)^2\rangle$ is abelian, but not finitely generated.

\subsubsection{Groups of the curves (2)-(2a)}
These curves are constructed as the 
curves (1)-(1a) with the following difference: 
One performs elementary transformations at the points 
$E_i\cap Q_i$ for $i=1,2,\dots, m$ 
followed by elementary transformations at the points 
$E_i\cap P_i$ for $i=m+1,m+2,\dots, n$. 
Finally, for $i=n,n+1\dots, 2n$ one applies elementary transformations
at some points $s_i\in Q_i\moins E_i$. 
The curves (2a) are obtained by setting $n=m$ in the above procedure.
\par
The same reasoning as in the case of the curves (1)-(1a) shows that
$\mu(Q_{m+1}):=\beta(\alpha\beta)^m$ is a meridian of 
$\widetilde{Q}$, and $\mu(P_{m+1}):=\alpha(\alpha\beta)^{m-n}$ 
is a meridian of $\widetilde{P}$. 
Setting $\mu(Q_{m+1})=\mu(P_{m+1})=1$, we obtain the presentation
$$
G_{(2)}:=\fg{\widetilde{C}}=
$$
$$
\langle \alpha, \beta, y\,|\,
\alpha=y^{-k}\beta^{k-1},\quad
[\beta,y^k]=\beta(\alpha\beta)^{m}=\alpha(\alpha\beta)^{n-m}=1
\rangle.
$$
Obviously, $\alpha\beta$ and $(\alpha\beta)\mon$ are singular
meridians at $\tilde{q}$ and $\tilde{p}$, respectively. A meridian of
$\widetilde{C}$ can be given as $a=\beta\mon y$. 
\par
For the curves (2a) we obtain,
$$
G_{(2a)}:=\fg{\widetilde{C}}=
$$
$$
\langle \alpha, \beta, y\,|\,
\alpha=y^{-k}\beta^{k-1},\quad
[\beta,y^k]=\alpha^m\beta^{m+1}=\alpha=1
\rangle.
$$
$$
=\langle \beta, y\,|\,
y^{-k}=\beta^{k-1},\quad
\beta^{n+1}=1
\rangle.
$$
Meridians in this case can be obtained from those of the case (2) by
putting $\alpha=1$. 
\par
Again, the element $y^k$ is central in the group $G_2$, and one has
$$
G_{(2)}(y^k)=
\langle y,\beta\,|\,
y^k=1,\quad\beta^j=1
\rangle=\integers_k*\integers_j,
$$
where this time
$j:=\mbox{g.c.d.}(1+mk, nk-mk+k-1)=\mbox{g.c.d.}(1+mk, nk+k)=
\mbox{g.c.d.}(1+mk, n+1)$.
\par
As for the curves (2a) we obtain,
$$
G_{(2a)}(y^k)=
\langle y,\beta\,|\,
y^k=1,\quad\beta^j=1
\rangle=\integers_k*\integers_j,
$$
where $j=\mbox{g.c.d.}(k-1, n+1)$.

\subsection{Groups of the curves (3)}
Let $C$ be the curve defined by the equation 
$xy^{k}-z^{k+1}=0$, the point $p$ be its cusp, $q$ its inflection
point, and $R$ be the line $\overline{pq}=\{z=0\}$. 
Let, as in the case (1), $P$, $Q$ be the tangent lines at $p$, $q$,
respectively.
\\ 
Blowing-up the inflection point $q=[0:1:0]$ 
we get the Hirzebruch surface $X_1$. 
For $i=1,2,\dots, n$, we perform elementary transformations 
at the points $E_i\cap R_i$, followed by elementary transformations
at some points $s_i\in Q_i\moins E_i$, and we end up with a Hirzebruch
surface $X_{2n+1}$ with $E_{2n+1}^2=1$. 
Let $\widetilde{C}$, $\widetilde{R}$, $\widetilde{Q}$ 
be the images of $C_{2n+1}$ ,$R_{2n+1}$, $Q_{2n+1}$ in 
$\projective$ under the contraction of $E_{2n+1}$.
Then $\widetilde{C}$ is a curve from the family (3).
One has the biholomorphic map
$$
\projective\moins(C\cup R\cup Q)
\stackrel{\simeq}{\longrightarrow}
\projective\moins(\widetilde{C}\cup \widetilde{R}\cup \widetilde{Q}),
$$
inducing an isomorphism 
$$
\fg{(C\cup R\cup Q)}\simeq
\fg{(\widetilde{C}\cup \widetilde{R}\cup \widetilde{Q})}.
$$
In order to find $\fg{(C\cup R\cup Q)}$, 
we shall use the projection
from the point $O:=P\cap Q$, as before. 
As the only points of intersection $R\cap C$ are $p$ and $q$, this
projection has only two singular fibers, namely $P$ and $Q$. 
Choose the base $B$ and the
generic fiber $F$ as in \fffigure{C_k1}. 
\figur{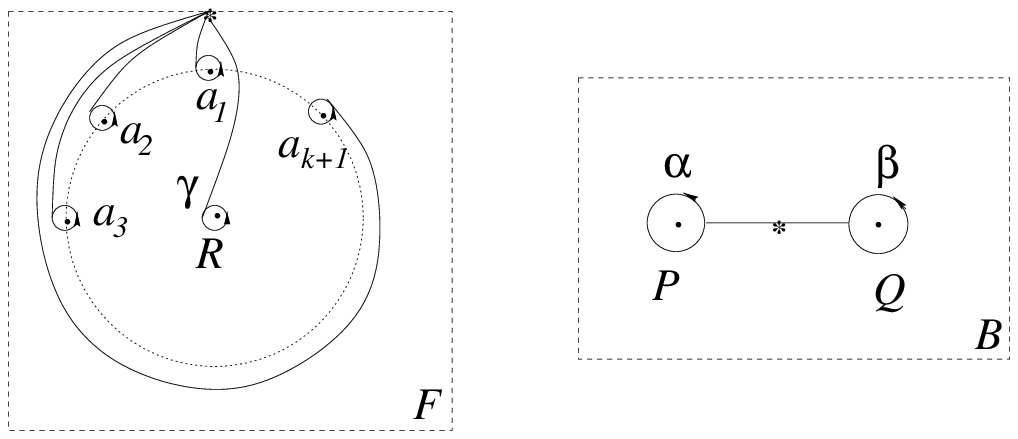}{}
Let $F^\prime:=F\moins(C\cup \{O\}\cup R)$, and choose the
generators $a_1,\dots a_{k+1}, \gamma$ 
for $\pi_1(F^\prime)$ as in \fffigure{C_k+1}. 
The monodromy relations around $Q$ are given by
\begin{center}
$\beta\mon a_i \beta =
\left\{
\begin{array}{ll}
      a_{i+1},& 1\leq i\leq k,\\
      (\delta\gamma)\;a_1\; 
      (\delta\gamma)^{-1}, & i=k+1,
\end{array}
\right.$\\
$\beta\mon \gamma \beta=$ $\qquad a_1\gamma a_1\mon.\qquad\qquad$\\
\end{center}
where $\delta:=a_{k+1} a_{k}... a_2 a_1$.
Observe that, as $\beta$ is a meridian of $Q$, 
and as the elementary transformations are applied 
at points $s_i\in Q_i\moins E_i$, $\beta$ is a meridian of $Q_i$, 
and thus it is a meridian of $\widetilde{Q}$. 
Imposing the relation $\beta=1$ in the relations
found above, we find that $a_1=a_2=\dots =a_{k+1}$, and $[a_1,\gamma]=1$.
Since the group $\fg{\widetilde{C}}$ is generated by these elements we
conclude that it is abelian.

\subsection{Groups of the curves (4)}
We begin by the curve $C:=xy^k-z^{k+1}$. 
Let $Q:=\{x=0\}$ be the tangent to $C$ at its flex $q:=[0:1:0]$,
and let $P$ be a line intersecting $C$ transversally at
its cusp $p:=[1:0:0]$, and such that $q\notin P$. 
Then by Bezout's theorem, $P$ intersects 
$C$ at one further point $r\neq q$. 
Let $O:=P\cap Q$. 
Blowing-up $O$, we get the Hirzebruch surface $X_1$,
with the horizontal section $E_1$ with $E_1^2=-1$. 
For $i=1,2,\dots, n$, apply elementary transformations 
at the points $r_i$ followed by elementary transformations
applied at the points $E_i\cap Q_i$. 
Then $E_{2n+1}^2=-1$; contracting it, we turn  back to the 
projective plane $\projective$. 
Let $\widetilde{C}$, $\widetilde{P}$, $\widetilde{Q}$ be the 
images of $C_{2n+1}$, $P_{2n+1}$ $Q_{2n+1}$ under this 
contraction. 
Then $\widetilde{C}$ is a curve of the family (4), and one has the
biholomorphism
$$
\projective\moins(C\cup P\cup Q)
\stackrel{\simeq}{\longrightarrow}
\projective\moins(\widetilde{C}\cup \widetilde{P}\cup \widetilde{Q}),
$$ 
inducing an isomorphism of the fundamental groups.
\\
To find the  group of $C\cup P\cup Q$, we shall use the
projection from the point $O$. 
In addition to $P$ and $Q$, this projection has a third singular fiber
$R$, which is a simple tangent to $C$ 
at a unique, smooth point of $C$. 
That $P$, $Q$, $R$ are the only singular fibers 
can be seen by looking at the dual picture.
Indeed, by the class formula, one has
$$
d^*=2(g-1+k+1)-(k-1)=k+1
$$
where $g=0$ is the genus of $C$, and $d^*$ is the degree of the dual
curve $C^*$. Now, the point $Q^*\in {\projective}^*$ is a cusp 
of $C^*$ with multiplicity $k$. 
Hence, the line $O^*$ which passes
through $Q^*$ should intersect $C^*$ transversally at a unique 
further point $R^*$, which is the dual of the simple tangent 
line $R$ from $O$ to $C$. 
\par
Now we apply the change of coordinates
$$
[x:y:z]\Rightarrow [x+y: y:z].
$$
In the new coordinates, the equation of $C$ reads as
$(x-y)y^k-z^{k+1}=0$. 
Let $L_\infty$ be the line $x=0$, and pass to the 
affine coordinates $(y/x,z/x)$ 
in $\complex^2=\projective\moins L_\infty$.
Recall that we have the freedom to choose $O$ 
(or, equivalently, $P$). 
So let $O=(1,z_0)$, where $z_0$ is a big real number.
\figur{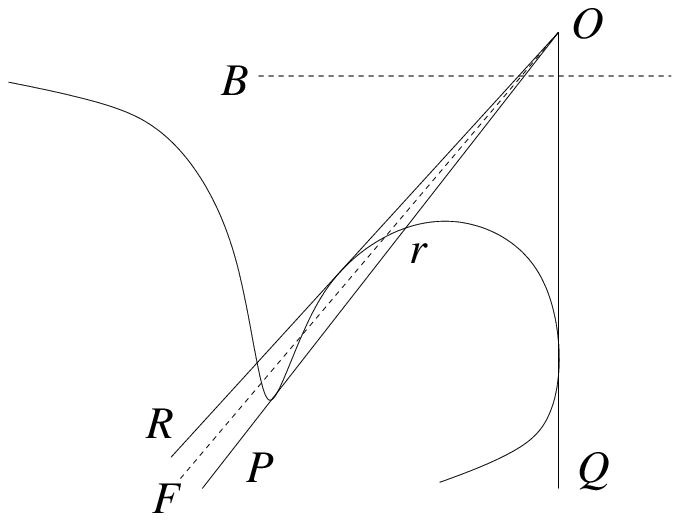}{}
The real picture of the configuration
$C\cup P\cup Q\cup R$ is shown in \fffigure{fenske4}. 
Choose the base fiber $F$, and the base of the projection 
as in \fffigure{fenske4}. 
Put $F^\prime:=F\moins (C\cup O)$ and
$B^\prime:=B\moins (P\cup Q\cup R)$. 
Take the generators $b,a_1,a_2, \dots, a_k$
for $\pi_1(F^\prime)$ as in \fffigure{fenske42}-I  and the 
generators $\alpha$, $\beta$,  $\gamma$ for $\pi_1(B^\prime)$ as 
in \fffigure{fenske42}-III.
\figur{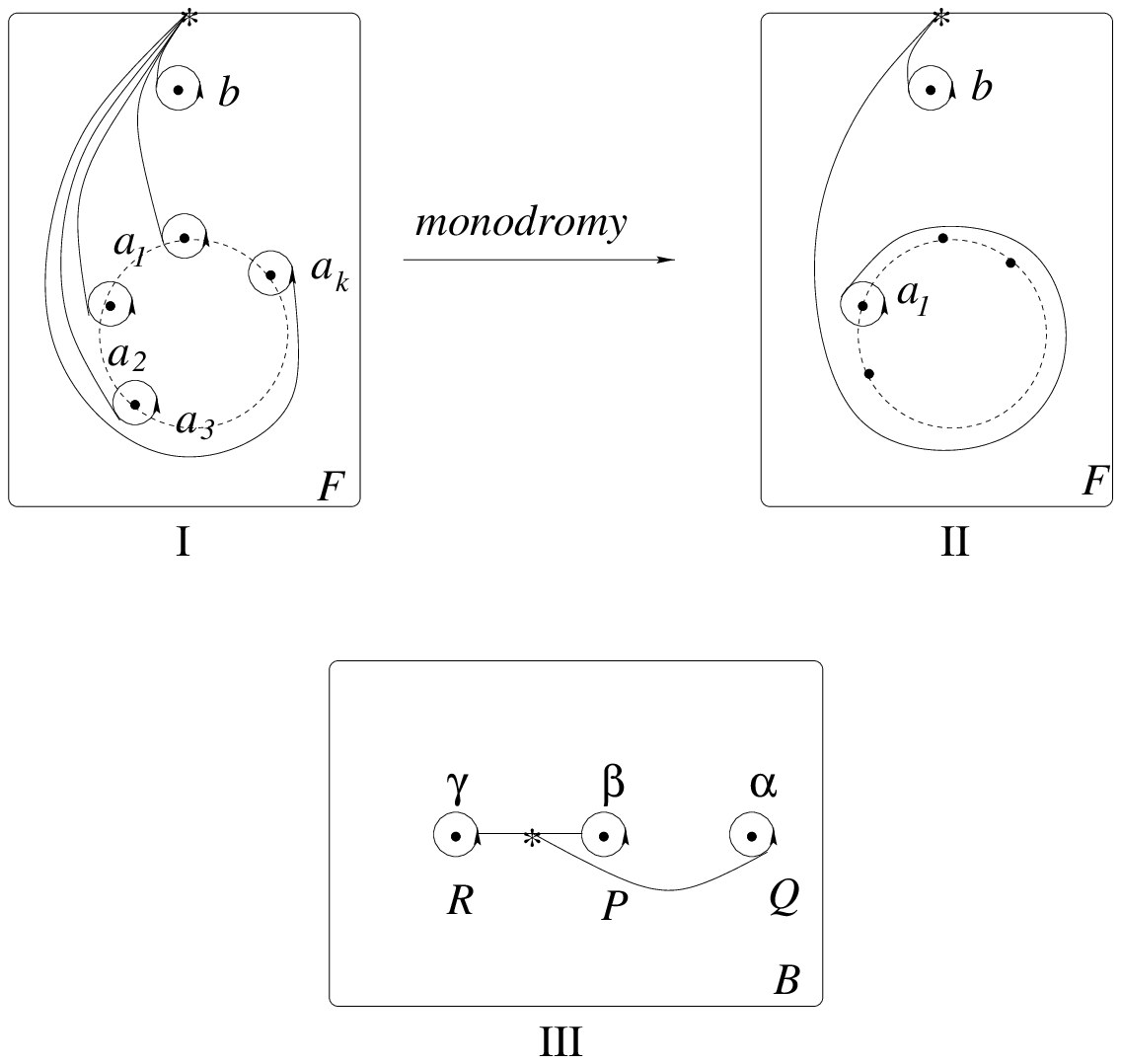}{}  
\par
The monodromy around $R$ yields (after setting $\gamma=1$),
$$
a_1=b\quad ({\cal R}_1),
$$
and the monodromy around $P$ gives
\begin{center}
$\beta\mon a_i \beta =
\left\{
\begin{array}{rl}
      \delta a_{i+1} \delta\mon \qquad\qquad& 1\leq i<k,\\
      \delta^2 a_{1} \delta\mtw  \qquad\qquad& i=k,
\end{array}
\right.
\quad ({\cal R}_2)$,\\
\end{center}
where $\delta:=a_ka_{k-1}\cdots a_1$. 
Hence, we have the presentation
$$
\fg{(C\cup P\cup Q\cup R)}=
\langle \beta, b, a_1,a_2,\dots ,a_k\,|\, ({\cal R}_1),\,({\cal R}_2)
\rangle
$$
Note that $\alpha\beta \delta b=1$, and that $[\alpha,\beta]=1$.
\par
An obvious application of Fujita's Lemma shows that 
$\beta b^n$ is a meridian of $\widetilde{P}$ and $\alpha(\alpha\beta)^n$
is a meridian of $\widetilde{Q}$. Therefore,
$$
G_{(4)}:=\fg{\widetilde{C}}=\langle \beta, b, a_1,a_2,\dots, a_k\,|\,
\beta b^n=\alpha(\alpha\beta)^n=1, \,({\cal R}_1),\,({\cal R}_2)
\rangle.
$$

\subsubsection{Study of the group}
By using $({\cal R}_2)$, one can express the generators 
$a_1, a_2,\dots ,a_k$ in terms of $b$ as follows
$$
a_i=(\beta\delta)^{-i+1}b (\beta\delta)^{i-1}\qquad 1\leq i\leq k.
$$
Then, the last relation in $({\cal R}_2)$ reads
$$
(\beta\delta)^{-k}b(\beta\delta)^k=\delta b\delta\mon\quad
({\cal R}_3).
$$
For $\delta$, one has
$$
\delta=a_ka_{k-1}\cdots a_1=(\beta\delta)^{-k}(\beta\delta b)^k
\quad ({\cal R}_4).
$$
Since $\alpha=(\beta\delta b)^{-k}$, the relation
$\alpha^{n+1}\beta^n=1$ becomes
$$
(\beta\delta b)^{n+1} b^{n^2}=1\quad ({\cal R}_5),
$$ 
where we have used $\beta=b^{-n}$.
This gives the presentation
$$
G_{(4)}=\langle \beta,\delta,b\,|\, 
\delta=(\beta\delta)^{-k}(\beta\delta b)^k,\quad
\beta=b^{-n},\quad
(\beta\delta b)^{n+1} b^{n^2}=1,\quad
(\beta\delta)^{-k}b(\beta\delta)^k=\delta b\delta\mon\rangle.
$$
Now put $x:=\beta\delta$ and $y:=b$. Then $\delta=\beta\mon x=y^n x$,
and one can rewrite the above presentation as
$$
G_{(4)}=\langle x,y\,|\, 
(xy)^{n+1}y^{n^2}=[y, x^ky^nx]=1,\quad x^ky^nx=(xy)^k\rangle,
$$
since $({\cal R}_3)$ becomes 
$$
x^{-k}yx^k=y^n x y x\mon y^{-n}
\,\Leftrightarrow\, [y, x^ky^nx]=1,
$$
and for $({\cal R}_4)$ one has
$$
y^n x=x^{-k}(xy)^k
\,\Leftrightarrow \, x^ky^nx=(xy)^k.
$$
Finally, $({\cal R}_5)$ is written as $(xy)^{n+1} y^{n^2}=1$.
Note that $y=b$ is a meridian of $\widetilde{C}$. 
\par
Obviously, the latter presentation is  equivalent to the presentation
$$
G_{(4)}=\langle x,y\,|\, 
(xy)^{n+1}y^{n^2}=[y, (xy)^k]=1,\quad x^ky^nx=(xy)^k\rangle \qquad (*)
$$
To simplify this presentation further, put $z:=xy$. 
Then $x=zy\mon$, and one obtains the presentation
$$
G_{(4)}=\langle z,y\,|\, z^{n+1}y^{n^2}=[y, z^k]=1,\quad
(zy\mon)^k y^n zy\mon=z^k \rangle
$$
It is readily seen from this presentation that the element $z^k$ is
central. Passing to the quotient by $z^k$ gives
$$
G_{(4)}(z^k)=\langle z,y\,|\,z^{n+1}y^{n^2}=1, \quad
(zy\mon)^{k+1}y^n=z^k=1 \rangle.
$$
Now put $j:=gcd(n+1,k)$. Then one has
$$
G_{(4)}(z^k)(y^n)=\langle z,y\,| z^j=y^n=(zy\mon)^{k+1}=1 \rangle=
T_{j,n,k+1},
$$
so that this latter group is big if 
\begin{equation}
\label{harmonic}
\frac{1}{j}+\frac{1}{n}+\frac{1}{k+1}< 1.
\end{equation}
\par
Obviously, $G_{(4)}(z^k)$, and hence also $G_{(4)}$ is abelian if 
$j=1$ or $n=1$. So suppose that $n,j\geq 2$. 
\par
First we consider the case $k=2$. 
Since $j\geq 2$, this forces $n\geq 3$ to be odd. 
If $n\geq 7$, then \ref{harmonic} is satisfied. 
In case $n=3$ or $n=5$, \ref{harmonic} is violated.
\par
Now assume $k=3$. Then $j=3$ by the assumption $j\geq 2$, which forces
$n+1\geq 3$ to be a multiple of $3$. For $n+1\geq 6$, 
\ref{harmonic} is not violated. For $n=2$, \ref{harmonic}
is violated. 
\par
For $k\geq 4$, first assume that $k$ is even. Then the least value 
that $j$ can take is $2$, and the least value that $n$ can take 
is $3$. If $n=3$, then $j=4$, and \ref{harmonic} is not violated.
But if $n\geq 4$, then \ref{harmonic} is not violated neither.
\par
If $k\geq 4$ is odd, then the least divisor of $k$ is $3$,  hence 
$j\geq 3$, and $k\geq 6$. In this case \ref{harmonic} is never
violated.
\par
This leaves the cases $(d,n,k)=(9,3,2)$, $(d,n,k)=(13,5,2)$, and
$(n,k)=(10,2,3)$ open.  Calculations with Maple show that these are finite
groups of order 72, 1560 and 240 respectively. 
\par
Notice that when $k=2$, the degree of the curves (4) is $2n+3$, so it
is interesting to compare their groups with the groups in Theorem
3.2.2. For $k=2$, the relation $x^2y^nx=(xy)^k$ in the presentation
(*) above becomes 
$$
x^2y^nx=xyxy\quad \Leftrightarrow xy^nx=yxy \quad (**).
$$
Now put 
$$
\alpha:=y, \quad \beta:=xy^{n-1}.
$$
Then 
$$
y=\alpha, \quad x=\beta\alpha^{1-n},
$$
and the relation (**) becomes the braid relation
$\beta\alpha\beta=\alpha\beta\alpha$.
On the other hand, the relation $[y, (xy)^k]=1$ in the presentation
(*) is written, in terms of $\alpha$, $\beta$, as $[\alpha,
\beta^n]=1$. Finally, the relation $(xy)^{n+1}y^{n^2}=1$ in (*) becomes
$$
(\beta\alpha^{2-n})^{n+1}\alpha^{n^2}=1.
$$
Hence, the presentation
$$
\langle \alpha, \beta \,|\, \beta\alpha\beta=\alpha\beta\alpha, 
\quad [\alpha,\beta^n]=(\beta\alpha^{2-n})^{n+1}\alpha^{n^2}=1\rangle.
$$
\subsection{Groups of the curves (5)}
Let $C$ be the curve $\{xy^{k-1}-z^k-z^{k-1}y\}$.
Its  unique 
singularity is a cusp at the point $p:=[1:0:0]$, and it has a
flex of order $k-1$ at the point $q:=[0:1:0]$. 
Let $P$, $Q$ be the tangents to $C$ at $p$ and $q$.
By Bezout's theorem, $Q$ intersect $C$ at a third point $r$.
Blowing-up the point $O:=P\cap Q$, we get the 
Hirzebruch surface $X_1$. 
Perform $n$ elementary transformations at 
$E_i\cap P_i$, followed by $n$ elementary 
transformations at the points $r_i$. 
One obtains the Hirzebruch surface $X_{2n+1}$
with $E_{2n+1}^2=-1$. Contraction of $E_{2n+1}$
gives the projective plane $\projective$; denote
by $\widetilde{C}$, $\widetilde{P}$, $\widetilde{Q}$
the images of $C$, $P$ and $Q$. 
Then $\widetilde{C}$ is a curve of the family (5).
\par
To find $\fg{(C\cup P\cup Q)}$, we shall use the projection 
from the point $O$. Evidently, $P$ and $Q$ are singular fibers 
of this projection. There is one further singular fiber say $R$, which 
is tangent to $C$ at a unique smooth point of $C$. Indeed, 
by the class formula, one has $d^*=k+1$. 
The point $Q^*\in{\projective}^*$ is a cusp of multiplicity
$k-1$. The line $O^*$ intersect $C^*$ at $P^*$, which is a smooth
point of $C^*$, and at the cusp $Q^*$.
By Bezout's theorem, $O^*$ should intersect $C^*$ at a third point
$R^*$ transversally, which is the point dual to the line $R^*$.
This reasoning shows also that there are no other singular fibers.
\par
In the affine coordinates $(x,z)$, the equation of $C$ reads
$x=z^k+z^{k+1}$, and it is easy to see that the third singular fiber
is the tangent line at $z=-k/(k+1)$. 
\figur{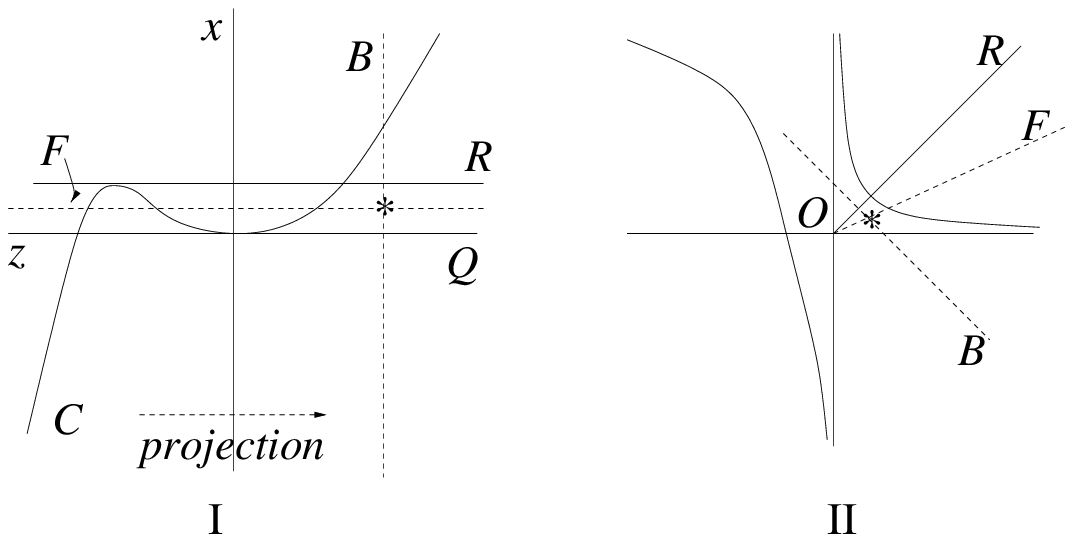}{}
For $k$ even, the situation is pictured in \fffigure{fenske5}.
Choose the base $B$ and the fiber $F$ as in the \fffigure{fenske5},
define $F^\prime$, $B^\prime$ as usual, and choose the generators for
their fundamental groups as in \fffigure{fenske52}.
\par
Set $\delta:=a_k a_{k-1}... a_2 a_1$. Then 
the monodromy around $Q$ yields the relations
$$
\beta^{-1} a_i \beta =
\left\{
\begin{array}{ll}
a_{i+1},& 1\leq i\leq k-1,\\
      \delta\;a_1\; \delta\mon, & i=k,
\end{array}
\right.
\quad[\beta,b]=1
$$
and the monodromy around $R$ gives
$$
a_k=b
$$
\figur{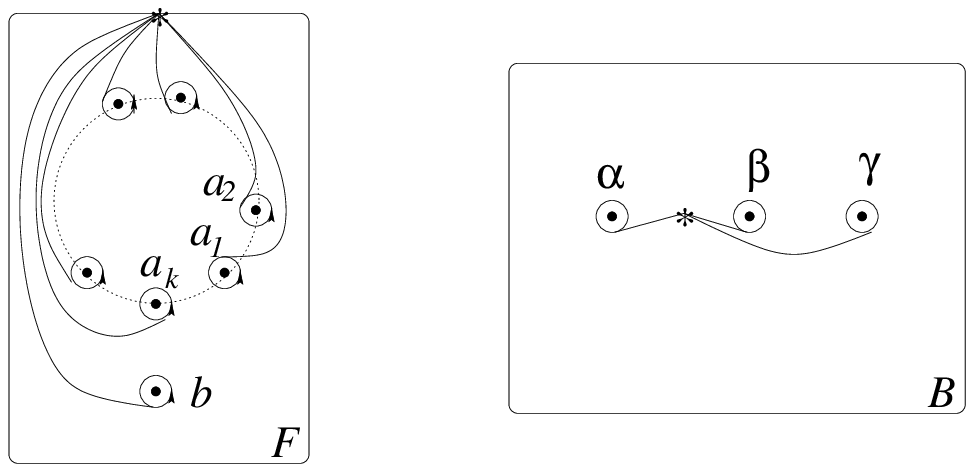}{}
\noindent
Now, the relation $\beta\mon a_{k-1}\beta=a_k=b$ implies $a_{k-1}=b$, since 
$[\beta, a_k]=1$. Similarly, one obtains $a_1=a_2=\cdots= a_k=b$.  
This shows that $\fg{(C\cup P\cup Q)}$ is abelian, which implies that 
$\fg{\widetilde{C}}$ is abelian, too. 
\subsection{Groups of the curves (6) in Theorem~\ref{table}}
Let $C$ be the curve $\{xy^{k+1}-z^{k+2}-z^{k+1}y=0\}$, and let 
$Q$ be its tangent line at the inflection point $q$. 
Denote by $S$ the line $\overline{pq}$.
The line $Q$ intersects $C$ at a second point $r$.
Blowing-up the point $q$, we get the Hirzebruch surface 
$X_1$ with the horizontal section $E_1$.
For $i=1,2,\dots, n$, perform elementary transformations at 
$E_i\cap S_i$, followed by elementary transformations at the
points $r_i$ for $i=n+1, n+2,\dots, 2n$. 
We end up with the Hirzebruch surface $X_{2n+1}$, with 
$E_{2n+1}=-1$. Contraction of $E_{2n+1}$ gives the
projective plane $\projective$, denote by 
$\widetilde{S}$, $\widetilde{Q}$, $\widetilde{C}$
the images of $S$, $Q$, and $C$ under this contraction.
\par
The calculation of $\fg{(C\cup S\cup Q)}$   
will be realized by using the projection from the 
point $O:=P\cap Q$, where $P$ is the tangent to
$C$ at its cusp as usual. The real picture is as
in \fffigure{fenske5}, except that we should 
take the line $S=\overline{pq}=\{z=0\}$ into 
consideration. Now, the only points of intersection 
of the line $S$ with $C$ are the points $p$ and $q$. 
Hence, the only singular lines of the projection from the point 
$O$ are $P$, $Q$, and $R$, where $R$ is the simple 
tangent line from $O$ to $C$, as shown in the previous section. 
\par
Choose a fiber $F$ and a base $B$ as in \fffigure{fenske5}-I,
and put $F^\prime:=F\moins(C\cup S\cup \{O\})$, 
$B^\prime:=B\moins(P\cup Q\cup R)$. The situation is illustrated in 
\fffigure{fenske6}.
\figur{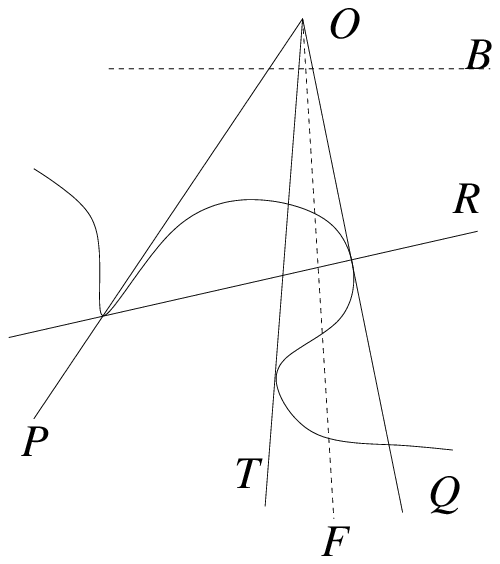}{}
Choose the generators $a_1, a_2, \dots ,a_{k+1}, b, \gamma$
for $\pi_1(F^\prime)$, and the generators 
$\alpha, \beta$ for $\pi_1(B^\prime)$ as in 
\fffigure{fenske62}.
\figur{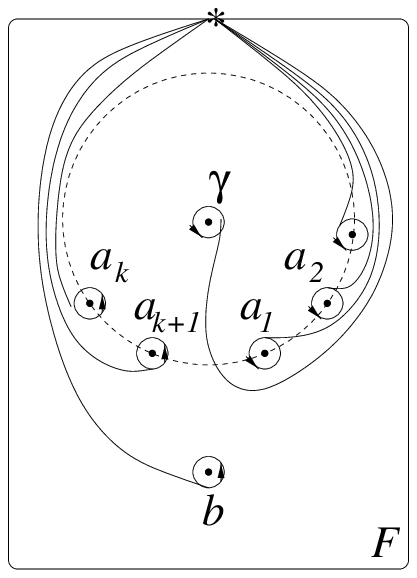}{}
\par
The monodromy around the singular fiber $Q$ gives
$$
\beta^{-1} a_i \beta =
\left\{
\begin{array}{rl}
a_{i+1},\qquad & 1\leq i\leq k,\\
(\delta\gamma) a_1 (\delta\gamma)^{-1}, & i=k+1,
\end{array}
\right.
\quad[\beta,b]=1,
$$
where $\delta:=a_{k+1}a_k\cdots a_1$. 
From the monodromy around $R$ we get,
$$
a_{k+1}=b.
$$ 
These relations yields  a presentation of the group  
$\fg{(C\cup S\cup  Q)}$.
\par
 Applying Fujita's Lemma to the elementary transformations applied
at the points $r_i$, we obtain the relation $\beta b^n=1$. 
Imposing this relation on the relations found above, we find that
$a_1=a_2=\cdots =a_{k+1}=b$, and $[b,\gamma]=1$. 
We conclude that the group is abelian.

\subsection{Groups of the curves (7) in Theorem~\ref{table}}
Let $C$ be the curve $xy^{k-1}-z^k=0$, with $p:=[1:0:0]$
as its cusp, $q:=[0:1:0]$ as its inflection point, and 
$P$, $Q$ as the tangents to $C$ at these points. 
Blow-up the point $O:=P\cap Q$, to get the Hirzebruch surface 
$X_1$, with the horizontal section $E_1$. 
Perform an elementary transformation at $q_1$, 
followed by an elementary transformation at $E_2\cap P_2$. 
Now for $i=3,4, \dots, n+2$ perform elementary transformations
applied at the points $E_i\cap P_i$, followed by 
elementary transformations applied at some points
$s_i\in Q_i\moins E_i$ for $i=n+3, n+4, \dots, 2n+2$.
 We end up with a Hirzebruch surface $X_{2n+3}$
with $E_{2n+3}=-1$. Contraction of $E_{2n+3}$ gives 
the projective plane $\projective$. 
\par
We shall use the projection from the point $O$ to 
find the fundamental group. This projection has, in addition to 
$P$ and $Q$, a third singular fiber $R$, which is a simple tangent
from $O$ to $C$ at a unique point.
The setting is same as in  Case 4, see \fffigure{fenske4}. 
Choose a base fiber close to $Q$, and take the generators for the base
and the fiber as in \fffigure{fenske42}. 
The monodromy around $Q$ gives the relations  
$$
\beta\mon a_i \beta =
\left\{
\begin{array}{ll}
a_{i+1},& 1\leq i\leq k-1,\\
\delta a_1 \delta\mon, & i=k,
\end{array}
\right.
$$
where $\delta:=a_ka_{k-1}\cdots a_1$.
\par
Now, without finding the monodromy around $R$ or $P$, notice that
an obvious application of Fujita's lemma gives $a_1\beta$ 
as a meridian of $Q_1$. Subsequent elementary transformations on 
$Q_i$ are applied at some points $s_i\in Q_i\moins E_i$, so 
that $a_1\beta$ stays to be a meridian of $Q_i$. 
Imposing the relation $a_1\beta=1$ on the above relations gives 
$\beta\mon =a_1=a_2=\cdots =a_k$. 
But, the group  $\fg{(C\cup P\cup Q)}$ is generated by 
the elements $a_1, a_2,\dots ,a_k, \beta$. 
Hence, the fundamental groups of the curves
of the family (7) are abelian.

\subsection{Groups of the curves (8) in Theorem~\ref{table}}
Let $C$ be the curve $xy^2-z^3=1$.
Pick an arbitrary smooth point $q\in C$ which is not the inflection
point of $C$, and let $Q$ be the tangent of $C$ at $q$.
Put $P:=\overline{pq}$, where $p$ is the cusp of $C$.
\par
By Bezout's theorem, the line $Q$ intersect 
$C$ at a second point, say $r$. 
\par
Blowing-up the point $q$, we get the Hirzebruch surface $X_1$, with 
the horizontal section $E_{1}$. 
For $i=1,2,\dots, n-1$, perform elementary transformations 
at the points $r_i$, 
followed by elementary transformations performed at the 
points $E_i\cap P_i$ for $i=n,n+1,\dots, 2n-2$.
We end up with the Hirzebruch surface $X_{2n-2}$ with 
$E_{2n-1}^2=-1$. Contraction of $E_{2n+1}$ gives the projective
plane $\projective$ back. 
Denote as usual by $\widetilde{P}$, $\widetilde{Q}$,
$\widetilde{C}$ the images of $P$, $Q$, $C$ under this 
contraction. Then $\widetilde{C}$ is a curve of the family (8).
\par
Let us show that the group $\fg{(C\cup P\cup Q)}$ is abelian, thereby 
showing that the groups of the curves of the family (8) are abelian.
\par
Consider the \textit{singular} projection from the cusp $p$. 
A generic fiber of this projection intersects $C\cup P\cup Q$ at 
two points, i.e.  $\fg{(C\cup P\cup Q)}$ is generated by two elements.
These two points meet at the point $r$, which is a transversal
intersection of $Q$ with $C$. This  implies that the corresponding 
generators commute. 
It follows that $\fg{(C\cup P\cup Q)}$ is abelian. $\Box$

\subsection{Groups of the curves in Theorem~\ref{fe5}}
Let $C$ be the curve $(yz-x^2)^2-x^3y=0$. 
\begin{lem}(Fenske~\cite{Fe2})
$C$ is a rational cuspidal quartic with cusps at the point
$p:=[0:0:1]$ of type $[2_2]$, and at the point
$r:=[0:1:0]$ of type $[2_1]$. 
This curve has an inflection point of order $3$ at the 
point $q:=[-576:-4096:135]$. 
\end{lem}

Let $P$ be the tangent to $C$ at the cusp $p$, and let $Q$ be the
inflectional tangent line at $q$. The cusp $p$ is the only
intersection point of $P$ with $C$, whereas $Q$ intersect 
$C$ at a second further point, say $s$. 
It is clear that the point $O:=P\cup Q$ does not lie on $C$.  
Blowing-up $O$, we get the Hirzebruch surface 
$X_1$, with a section $E_1$ such that $E_1^2=-1$. 
Now for $i=1,2,\dots ,n$ perform 
elementary transformations at the points 
$s_i$, followed by elementary transformations applied at the
points $E_i\cup P_i$ for $i=n+1,n+2,\dots, 2n+1$. 
We end up with the Hirzebruch surface $X_{2n+1}$ with 
$E_{2n+1}^2=-1$. Contraction of $E_{2n+1}^2=-1$ gives the
projective plane $\projective$.
Denote by $\widetilde{C}$, $\widetilde{P}$, $\widetilde{Q}$ the images
of $C_{2n+1}$, $P_{2n+1}$, $Q_{2n+1}$ under this contraction. 
Then $\widetilde{C}$ is the desired curve of degree $3n+4$.
 \par
To find the fundamental group of $C\cup P\cup Q$, 
 we shall use the projection from the center $O$. 
Let $R$ be the line $\overline{Or}$. Then, clearly 
$P$, $Q$ and $R$ are  singular fibers of this projection. 
That there are no other singular fibers can be seen by looking
at the dual picture. By the class formula, the degree of the curve 
$C^*$ dual to $C$ is $4$. 
The line $O^*$ intersects $C^*$ at its simple cusp $Q^*$ with
multiplicity $2$. Since the line $P$ intersects $C$ with multiplicity 
$4$ at the cusp $p$, $O^*$ intersects $C^*$ with multiplicity 
$2$ at the cusp $P^*$ of $C^*$. (That $P^*$ is a cusp of $C^*$ of
multiplicity $2$ can be seen by using the parameterization 
$[t^2:t^4: 1+t]$ of $C$.)
By Bezout's theorem, this accounts for all the
intersection points of the line $O^*$ with $C^*$.
\par
To get a better picture of the curve, we apply the transformation 
$[x:y:z:]\rightarrow [x:y:y+z]$, and then pass to the affine
coordinate system in $\complex^2=\projective\moins L_\infty$, where
$L_\infty =\{z=0\}$. 
In these coordinates, $C$ is parameterized as  
$(t^2/(1+t+t^4),t^4/(1+t+t^4))$, the cusp 
$p$ is the point $(0,0)$, the cusp $r$ is the point $(0,1)$, and the
flex $q$ is the point $(576/3961, 4096/3961)$. 
It turns out that the point $s$, the second point of intersection of
$Q$ with $C$, is real. 
The configuration $C\cup P\cup Q$ is pictured
in \fffigure{fensked4}. 
\figur{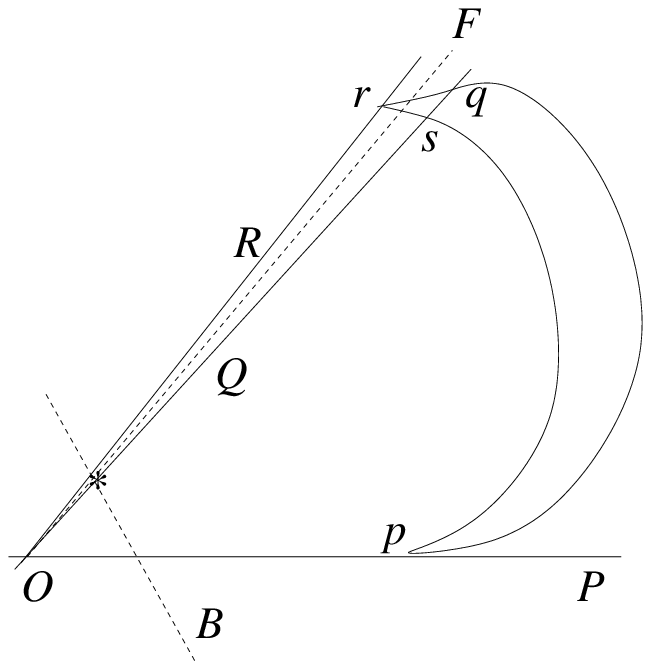}{}
Pick $F$, $B$ as in \fffigure{fensked4},
put $F^\prime:=F\moins (C\cup \{O\})$, 
$B^\prime:=B\moins (P\cup Q \cup R)$, 
and choose the generators $a,b,c,d$ of $F^\prime$ 
and the generators $\alpha, \beta, \gamma$ of $B^\prime$ as in 
\fffigure{fenskeson}. 
\figur{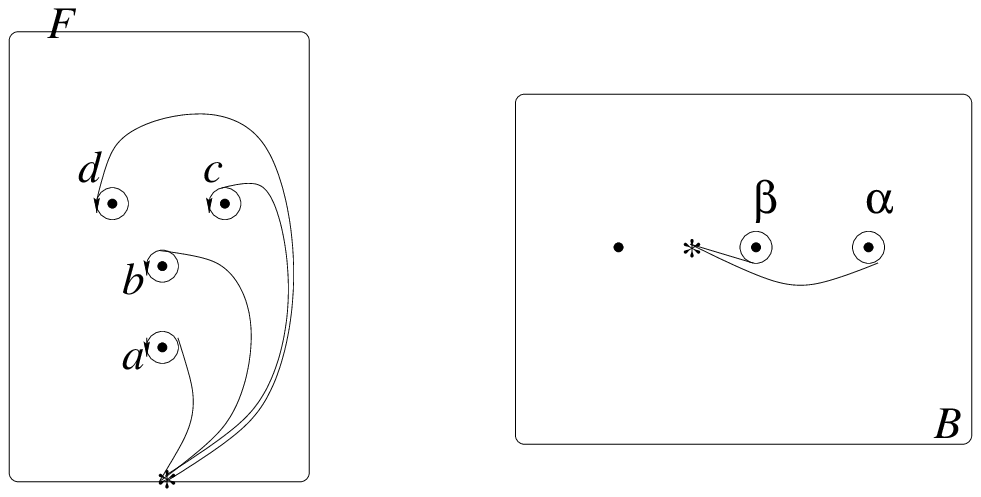}{}
\par
The monodromy around $Q$ gives the relations
$$
\beta\mon b \beta =c,\quad \beta\mon c\beta =d, 
\quad \beta\mon d \beta= dcb(dc)\mon, 
$$
and the relation 
$$
\beta\mon a \beta =a \Leftrightarrow [\beta, a]=1.
$$
Using these relations, 
one can express the generators $c$, $d$ in terms of $b$
and $\beta$, and one easily deduces the relation
$$
(b\beta)^3=(\beta b)^3.
$$
Finally, the monodromy around $P$ gives the cusp relation
$$
aba=bab.
$$
This completes the presentation of $\fg{(C\cup P\cup Q)}$. 
Note that the loop $\alpha$ can be expressed by using the  relation
$$
dcba\alpha\beta=1=\beta\mtw (b\beta)^3 a\alpha,
$$
which implies that $[\alpha,\beta]=1$, this latter can also be derived
by an application of Fujita's lemma at the point  $O$. 
\par
By Fujita's lemma, $a^n \beta$ is a meridian of $\widetilde{Q}$, 
and $\alpha^{n+1}\beta^n$ is a meridian of $\widetilde{P}$. 
Imposing the corresponding relations on the above presentation we get
$$
\fg{\widetilde{C}}=
\langle a,b,\alpha, \beta \,|\quad
(b\beta)^3=(\beta b)^3 \quad ({\cal R}_1),\quad
aba=bab \quad ({\cal R}_2),\quad
$$
$$
\beta a^n=1\quad ({\cal R}_3),\quad 
 \alpha^{n+1}\beta^n=1\quad ({\cal R}_4),\quad
\beta\mtw (b\beta)^3 a\alpha=1 \quad ({\cal R}_5)\rangle
$$
Note that the relation $[\beta ,a]=1$ follows from the relation
$({\cal R}_3)$ and therefore is not written in the above presentation. 
\par
In order to simplify this presentation, put $b\beta=:\gamma$.
Then the relation $({\cal R}_1)$ can be expressed as 
$[\beta, \gamma^3]=1$, or, substituting $\beta=a^{-n}$, as 
$[a^n, \gamma^3]=1$. On the other hand, $b=\gamma\beta\mon=\gamma a^n$,
so that the cusp relation $({\cal R}_2)$ becomes
$$
a \gamma a^na=\gamma a^n a \gamma a^n\Leftrightarrow
a\gamma a= \gamma a^{n+1} \gamma.
$$
Using $({\cal R}_5)$, we can express $\alpha$ as follows:
$$
\alpha=a\mon \gamma\mth a^{-2n}.
$$
Substituting this in $({\cal R}_4)$ we get, using $[a^n, \gamma^3]=1$,
$$
(a\mon \gamma\mth a^{-2n})^{n+1} a^{n^2}=1
\Leftrightarrow a^{3n^2+2n}(\gamma^3 a)^{n+1}=1.
$$
So we have obtained the presentation
$$
\fg{\widetilde{C}}=
\langle a, \gamma\,|\, a\gamma a= \gamma a^{n+1} \gamma,\quad
[a^n,\gamma^3]=a^{3n^2+2n}(\gamma^3 a)^{n+1}=1 \rangle.
$$
Put $G:=\fg{\widetilde{C}}$. Then one has
$$
G(n)=G(a^n)=
\langle a, \gamma\,|\, a\gamma a=\gamma a \gamma,\quad 
(\gamma^3 a)^n=1\rangle.
$$
Applying the transformation $x:=a\gamma$, $y:=a\gamma a$, and imposing
the relation $x^3=1$ we get a surjection $G\onto T$, where
$$
T:=
\langle x,y\,|\, y^2=x^3=x^{2(n+1)}=(xy\mon)^n=1 \rangle,
$$
so that if $3|(n+1)$, then one has a surjection
$G(n)\onto T_{2,3,n}$
onto the triangle group, which implies that $G(n)$ is big for $n\geq 7$.
 \par
A more economical way of writing the last relation in the presentation
of  $G$ is as follows: one has
$a^{n+1}=\gamma\mon a\gamma a\gamma\mon$. Hence,
$$
a^{3n^2+2n}(\gamma^3a)^{n+1}=a^{2n^2+n-1}a(\gamma^3 a^{n+1})^3
$$
$$
=(a^{n+1})^{2n-1}a\gamma^2(\gamma a^{n+1}\gamma\cdot \gamma\cdot
\gamma a^{n+1}\cdots \gamma a^{n+1}\gamma)\gamma\mon
$$
$$
=(\gamma\mon a\gamma a\gamma\mon)^{2n-1}a\gamma^2 (a\gamma
a)^{n+1}\gamma\mon,
$$
so that one can replace the last relation by the relation
$$
(\gamma\mon a\gamma a\gamma\mon)^{2n-1}a\gamma (\gamma a)^{2n+2}\gamma\mon.
$$
\textbf{Remark.} By the following lemma, the group 
$G$ is actually a quotient of the braid group
$\braids_3$ on three strands.
\begin{lem}
For any $n\in \integers$, one has
$
\langle a, b\,|\, aba=ba^{n+1}b\rangle \simeq \braids_3.
$
\end{lem}
\proof
Applying the transformation $(a,b)\rightarrow (x:=a, y:=ba^n)$,
with inverse $(x,y)\rightarrow (a=x, b=yx^{-n})$, the above relation
becomes
$$
x\, yx^{-n}\, x\,=\, yx^{-n}\,x^{n+1}\,yx^{-n},
$$
which is nothing but the braid relation $xyx=yxy. \qquad \Box$

\end{document}